\documentclass[11pt]{amsart}

\usepackage[a4paper,margin=2.8cm]{geometry}
\usepackage{amsmath,amssymb,amsthm,mathtools}
\usepackage{enumitem}
\usepackage{microtype}
\usepackage{xcolor}
\usepackage{tikz}
\usetikzlibrary{arrows.meta,calc}
\usepackage{tikz-3dplot}
\usepackage{xurl}
\usepackage[colorlinks=true,linkcolor=blue!50!black,citecolor=green!35!black,urlcolor=blue!50!black]{hyperref}
\usepackage[capitalise,nameinlink]{cleveref}

\theoremstyle{plain}
\newtheorem{theorem}{Theorem}[section]
\newtheorem{proposition}[theorem]{Proposition}
\newtheorem{lemma}[theorem]{Lemma}
\newtheorem{corollary}[theorem]{Corollary}
\theoremstyle{definition}
\newtheorem{definition}[theorem]{Definition}
\newtheorem{example}[theorem]{Example}
\newtheorem{question}[theorem]{Problem}
\theoremstyle{remark}
\newtheorem{remark}[theorem]{Remark}

\newcommand{\R}{\mathbb{R}}
\newcommand{\N}{\mathbb{N}}
\newcommand{\Fix}{\operatorname{Fix}}
\newcommand{\diam}{\operatorname{diam}}
\newcommand{\bmeet}{\textstyle\bigwedge}
\newcommand{\bjoin}{\textstyle\bigvee}
\newcommand{\dotminus}{\mathbin{\dot{-}}}

\title[Smyth completeness and nonexpansive iteration]{Smyth completeness and nonexpansive iteration\\ in $q$-hyperconvex $T_0$-quasi-metric spaces}
\author{Collins Amburo Agyingi}
\address[C.~A.~Agyingi]{Department of Mathematical Sciences, University of South Africa (Unisa), Pretoria, South Africa}

\author{Ya\'e Ulrich Gaba}
\address[Y.~U.~Gaba]{AI Research and Innovation Nexus for Africa (AIRINA Labs), B\'enin; African Center for Advanced Studies (ACAS), Cameroon}

\thanks{Corresponding author: C.~A.~Agyingi.}

\date{\today}
\subjclass[2020]{Primary 47H10, 54E35; Secondary 06B23, 54C15, 18D20}
\keywords{$T_0$-quasi-metric, specialization order, $q$-hyperconvexity, Isbell-convexity, nonexpansive map, complete lattice, Kleene iteration}

\begin{document}

\begin{abstract}
We study the specialization order $\le_d$, defined by $x\le_d y$ if and only if $d(x,y)=0$, of a $q$-hyperconvex $T_0$-quasi-metric space $(X,d)$. We show that every nonempty bounded subset $A$ of $X$ has an infimum $\bigwedge A$ characterized by the metric formula $d(w,\bigwedge A)=\sup_{a\in A}d(w,a)$ for all $w\in X$, and dually a supremum, and that powers and copowers exist. For bounded $X$ this gives an elementary proof, with explicit metric formulas, of the known fact that Isbell-convex Lawvere metric spaces are complete in the sense of enriched category theory; the formulas themselves and the extension to bounded subsets of unbounded $q$-hyperconvex spaces are not covered by the enriched-categorical argument, which requires radii in $[0,\infty]$. Consequently no $T_1$-quasi-metric space with at least two points is $q$-hyperconvex, and the fixed points of every order-preserving self-map of a bounded $q$-hyperconvex space form a nonempty complete lattice. We show that the fixed point set of a nonexpansive self-map need not be a sublattice and that the specialization order need not be distributive. Our main result concerns the continuity of nonexpansive maps. Nonexpansive self-maps of bounded $q$-hyperconvex spaces need not preserve suprema of increasing sequences, and the supremum of the orbit of the least element need not be a fixed point. We prove that for a bounded $q$-hyperconvex $T_0$-quasi-metric space the following are equivalent: the space is Smyth complete; every increasing sequence converges to its supremum in the symmetrized metric; every nonexpansive self-map preserves suprema of increasing sequences; every nonexpansive self-map preserves suprema of directed sets; every nonexpansive self-map is Yoneda-continuous. In this case Kleene iteration from the least element converges to the least fixed point; this applies for instance when the symmetrized metric is compact.
\end{abstract}

\maketitle

\section{Introduction}\label{sec:intro}

K\"unzi and his school developed an asymmetric counterpart of the theory of hyperconvex metric spaces \cite{AP1956,Isbell1964,EK2001}: the theory of $q$-hyperconvex, or Isbell-convex, $T_0$-quasi-metric spaces and of the Isbell-hull \cite{KKO2012,Agyingi2014,Haihambo2013}. Fixed point theorems for nonexpansive maps followed \cite{KO2012,MO2013}. In these results the existence of fixed points is obtained by passing to the symmetrization $d^s=\max\{d,d^{-1}\}$: nonexpansiveness symmetrizes, and $(X,d^s)$ is hyperconvex when $(X,d)$ is $q$-hyperconvex \cite{KKO2012}, so the metric theorems of Sine, Soardi and Baillon apply. A survey of this theory, with complete proofs, is given in the companion paper \cite{Survey}.

In this article we exploit an ingredient that is invisible to $d^s$: the \emph{specialization order} $x\le_d y\iff d(x,y)=0$. Nonexpansive maps are order-preserving, and order intervals are balls of radius $0$. More is true: the specialization order is lost under every pseudometric $d'\ge d^s=\max\{d,d^{-1}\}$, since any such $d'$ satisfies $d'(x,y)\ge d(x,y)$ and $d'(x,y)\ge d(y,x)$, so $d'(x,y)=0$ forces $x=y$ by $T_0$. The specialization order of every such $d'$ is the identity relation, and the asymmetric content $\le_d$ disappears under every bona fide symmetrization. Fujii observed that Isbell-convexity is a geometric form of categorical completeness for Lawvere metric spaces \cite{Fujii2020}; for bounded spaces this implies that $(X,\le_d)$ is a complete lattice, so that for bounded $X$ \cref{thm:meets} below extracts the consequences of Fujii's theorem for the specialization order. What is new here: the proof is elementary and provides explicit metric formulas reused modularly in the proofs of \cref{thm:tarski-qh,prop:cne-qh,prop:comparable}, and the result extends to bounded subsets of unbounded $q$-hyperconvex spaces (\cref{thm:meets,prop:powers}), which lie outside the enriched-categorical argument because the latter requires radii in $[0,\infty]$. The lattice structure is then what drives the fixed-point results.

The specialization order has already been used to obtain fixed point theorems in quasi-metric spaces. Rutten proved that a continuous nonexpansive self-map $T$ of a complete quasi-metric space with $x_0\le_dTx_0$ has a fixed point above $x_0$ \cite{Rutten1996}. Valero and his coauthors developed Kleene-type fixed point theorems in quasi-metric spaces endowed with the specialization order, introduced order-completeness and order-continuity notions adapted to this setting, showed that Smyth completeness implies these completeness notions, and applied the results to the asymptotic complexity analysis of algorithms \cite{LRV2018,SV2018,EMV2019}. These works treat general quasi-metric spaces and impose order-continuity on the maps. Our point of view is different: in a $q$-hyperconvex space the specialization order is a (conditionally) complete lattice with metric lattice operations, so fixed points of all order-preserving maps exist without any continuity, and the setting lets us decide exactly when nonexpansive maps are order-continuous (\cref{thm:kleene}).

Our main results are the following. Let $(X,d)$ be a bounded $q$-hyperconvex $T_0$-quasi-metric space.
\begin{enumerate}[label=(\arabic*)]
\item The fixed points of every order-preserving self-map of $X$, continuous or not, form a nonempty complete lattice (\cref{thm:tarski-qh}). For nonexpansive maps this gives an order-theoretic proof of the existence part of the fixed point theorem of K\"unzi and Olela Otafudu \cite[Theorem~3.3]{KO2012}.
\item The fixed point set of a nonexpansive self-map need not be a sublattice of $X$, and $(X,\le_d)$ need not be distributive (\cref{thm:not-sublattice,prop:M3}).
\item Nonexpansive self-maps need not preserve suprema of increasing sequences, and the supremum of the orbit $(T^n\bot)$ of the least element need not be a fixed point (\cref{ex:shift}). The following are equivalent (\cref{thm:kleene}): $X$ is Smyth complete; every increasing sequence in $X$ converges to its supremum in $(X,d^s)$; every directed subset of $X$ converges to its supremum in $(X,d^s)$; every nonexpansive self-map of $X$ preserves suprema of increasing sequences, or of directed sets; every nonexpansive self-map of $X$ is Yoneda-continuous in the sense of Bonsangue, van Breugel and Rutten \cite{BvBR1998}. The proof uses copowers and the injectivity of $X$. Under these conditions $T^n\bot$ converges to the least fixed point of every nonexpansive $T$; they hold in particular when $(X,d^s)$ is compact (\cref{prop:compact}).
\end{enumerate}
We also show that a self-map of a $q$-hyperconvex space which is nonexpansive on comparable pairs is nonexpansive (\cref{prop:cne-qh}).

Throughout, $t\dotminus r=\max\{t-r,0\}$.

\section{Preliminaries}\label{sec:prelim}

This section fixes notation and recalls the minimum necessary to parse what follows. For a self-contained treatment of the underlying theory with complete proofs, test-case Isbell-hull computations, and the comparison of asymmetric completeness notions, we refer to the companion survey~\cite{Survey}.

\emph{Quasi-metrics.} A \emph{quasi-pseudometric} on a set $X$ is a map $d\colon X\times X\to[0,\infty)$ with $d(x,x)=0$ and $d(x,z)\le d(x,y)+d(y,z)$ \cite[Def.~2.1]{Survey}. It is a \emph{$T_0$-quasi-metric} if $d(x,y)=0=d(y,x)$ implies $x=y$, and a \emph{$T_1$-quasi-metric} if $d(x,y)=0$ implies $x=y$. The \emph{conjugate} is $d^{-1}(x,y)=d(y,x)$ and the \emph{symmetrization} $d^s=\max\{d,d^{-1}\}$ is a metric whenever $d$ is $T_0$; $(X,d)$ is \emph{bicomplete} if $(X,d^s)$ is complete. We write $C_d(x,r)=\{y:d(x,y)\le r\}$ and $C_{d^{-1}}(x,s)=\{y:d(y,x)\le s\}$. A map is \emph{nonexpansive} if $d(fx,fy)\le d(x,y)$ throughout. A set $A$ is \emph{bounded} if $\diam A=\sup\{d(a,a'):a,a'\in A\}<\infty$.

\emph{Specialization order.} The \emph{specialization order} of a $T_0$-quasi-metric space is $x\le_d y\iff d(x,y)=0$ \cite[Def.~2.7]{Survey}; it is a partial order, and $\le_{d^{-1}}$ is its reverse. We put ${\uparrow}a=C_d(a,0)$ and ${\downarrow}a=C_{d^{-1}}(a,0)$. A subset $D$ is \emph{directed} if it is nonempty and any two of its elements have an upper bound in $D$; we regard $D$ as a net indexed by itself. The following basic facts are proved in \cite[Lemma~2.8]{Survey}.

\begin{lemma}\label{lem:spec}
Let $(X,d)$ be a $T_0$-quasi-metric space.
\begin{enumerate}[label=(\alph*)]
\item Every nonexpansive map between $T_0$-quasi-metric spaces is order-preserving for $\le_d$.
\item For $a\in X$ the maps $d(a,\cdot)$ and $d(\cdot,a)$ are $1$-Lipschitz on $(X,d^s)$; hence ${\uparrow}a$ and ${\downarrow}a$ are closed in $(X,d^s)$.
\item $\le_d$ is the identity relation if and only if $d$ is a $T_1$-quasi-metric; in particular $\le_d$ is trivial whenever $d$ is a metric.
\end{enumerate}
\end{lemma}

\emph{Smyth completeness and Yoneda limits.} A net $(x_i)_{i\in I}$ is \emph{left $K$-Cauchy} if for every $\varepsilon>0$ there is $i_\varepsilon$ with $d(x_i,x_j)<\varepsilon$ whenever $i_\varepsilon\le i\le j$; the space is \emph{Smyth complete} if every left $K$-Cauchy net converges in $(X,d^s)$, equivalently, by \cite[Prop.~1]{RT2015}, if every left $K$-Cauchy sequence does \cite[Def.~3.1, Lem.~3.2]{Survey}. A point $y$ is a \emph{Yoneda limit} of a left $K$-Cauchy net $(x_i)$ if $d(y,w)=\inf_i\sup_{j\ge i}d(x_j,w)$ for all $w\in X$ \cite{BvBR1998}; it is unique when it exists. Nonexpansive maps send left $K$-Cauchy nets to left $K$-Cauchy nets, and a nonexpansive $T$ is \emph{Yoneda-continuous} if, whenever $y$ is a Yoneda limit of $(x_i)$, $Ty$ is a Yoneda limit of $(Tx_i)$ \cite{BvBR1998,GLN2017}. Convergence in $(X,d^s)$ implies Yoneda convergence, since $|d(x_j,w)-d(y,w)|\le d^s(x_j,y)$.

\emph{$q$-hyperconvexity.} Following Kemajou, K\"unzi and Olela Otafudu \cite{KKO2012} (see \cite[Def.~5.2]{Survey}), a quasi-pseudometric space $(X,d)$ is \emph{$q$-hyperconvex} (or \emph{Isbell-convex}) if, for every family $(x_i)_{i\in I}$ in $X$ and all families $(r_i)_{i\in I}$, $(s_i)_{i\in I}$ of nonnegative reals,
\begin{equation}\label{eq:qh}
d(x_i,x_j)\le r_i+s_j\ \text{ for all } i,j\in I\ \Longrightarrow\ \bigcap_{i\in I}\bigl(C_d(x_i,r_i)\cap C_{d^{-1}}(x_i,s_i)\bigr)\ne\emptyset .
\end{equation}
The hypothesis in \eqref{eq:qh} is necessary ($d(x_i,x_j)\le d(x_i,z)+d(z,x_j)\le r_i+s_j$ for $z$ in the intersection), and for $I=\emptyset$ the intersection is $X$, so $q$-hyperconvex spaces are nonempty. When $d$ is a metric and $r_i=s_i$, \eqref{eq:qh} is the hyperconvexity of Aronszajn and Panitchpakdi \cite{AP1956,EK2001}. The symmetrization $(\R,u)$ with $u(x,y)=x\dotminus y$, whose double balls are the compact intervals $[x-r,x+s]$ and whose specialization order is the usual order on $\R$, is the canonical unbounded $q$-hyperconvex example \cite[Ex.~2.11, Lem.~5.5]{Survey}. The basic structural facts we use are collected in the next proposition; the proofs are an exercise with \eqref{eq:qh} and are given in \cite[Props.~5.4, 5.7]{Survey}.

\begin{proposition}[{\cite[Prop.~2]{KKO2012}; \cite[Prop.~3.1]{KO2012}; \cite[Props.~5.4, 5.7]{Survey}}]\label{prop:basic}
Let $(X,d)$ be a $q$-hyperconvex quasi-pseudometric space.
\begin{enumerate}[label=(\alph*)]
\item $(X,d^{-1})$ is $q$-hyperconvex.
\item $(X,d^s)$ is a hyperconvex pseudometric space; if $d$ is $T_0$, it is a complete metric space.
\item If $(x_i,r_i,s_i)_{i\in I}$ satisfies the hypothesis of \eqref{eq:qh}, then $D=\bigcap_i\bigl(C_d(x_i,r_i)\cap C_{d^{-1}}(x_i,s_i)\bigr)$ is $q$-hyperconvex.
\end{enumerate}
\end{proposition}

\emph{Injective hull.} The categorical counterpart of $q$-hyperconvexity is injectivity in the category of $T_0$-quasi-metric spaces and nonexpansive maps.

\begin{theorem}[{Kemajou--K\"unzi--Olela Otafudu \cite[Thm.~1]{KKO2012}; \cite[Thm.~5.10, Cor.~5.11]{Survey}}]\label{thm:injective}
A $T_0$-quasi-metric space is $q$-hyperconvex if and only if it is injective in the category of $T_0$-quasi-metric spaces and nonexpansive maps. In particular, every $q$-hyperconvex subspace $Y$ of a $T_0$-quasi-metric space $X$ is the image of a nonexpansive retraction $X\to Y$.
\end{theorem}

Every $T_0$-quasi-metric space $X$ admits an \emph{Isbell-hull} $Q(X)$: the set of minimal function pairs $(f,g)$ with $f,g\colon X\to[0,\infty)$ and $g(x)+f(y)\ge d(x,y)$, under the asymmetric distance $N((f,g),(f',g'))=\max\{\sup_y(f(y)\dotminus f'(y)),\sup_x(g'(x)\dotminus g(x))\}$ \cite{KKO2012,Haihambo2013}; it is $q$-hyperconvex, and $X$ embeds isometrically via $x\mapsto(d(x,\cdot),d(\cdot,x))$. Minimality of $(f,g)$ amounts to $f(y)=\sup_x(d(x,y)\dotminus g(x))$ and $g(x)=\sup_y(d(x,y)\dotminus f(y))$. Explicit computations of Isbell-hulls of two- and three-point spaces are given in \cite[\S5]{Survey}; in particular, the Isbell-hull of the symmetric two-point space is isometric to $([0,1]^2,u_\infty)$ \cite[Ex.~5.6]{Survey}.

\emph{Products.} We use the following elementary fact repeatedly, for instance to establish that $([0,1]^K,u_\infty)$ is a bounded $q$-hyperconvex $T_0$-quasi-metric space whose specialization order is the coordinatewise order. The proof follows by applying \eqref{eq:qh} coordinatewise; see \cite[Ex.~2.12, Lem.~5.5]{Survey}.

\begin{lemma}\label{lem:products}
Let $(X_k,d_k)_{k\in K}$ be a family of $q$-hyperconvex quasi-pseudometric spaces such that $\sup_k\diam X_k<\infty$, and let $d(x,y)=\sup_kd_k(x_k,y_k)$ on $X=\prod_kX_k$. Then $(X,d)$ is $q$-hyperconvex. In particular $([0,1]^K,u_\infty)$ with $u_\infty(x,y)=\sup_k(x_k\dotminus y_k)$ is a bounded $q$-hyperconvex $T_0$-quasi-metric space, and $\le_{u_\infty}$ is the coordinatewise order.
\end{lemma}

\section{Metric meets, joins, powers and copowers}\label{sec:lattice}

Infima and suprema for the specialization order admit explicit metric formulas; powers and copowers are constructed, and the lattice-theoretic consequences used later are derived from them.

\begin{theorem}\label{thm:meets}
Let $(X,d)$ be a $q$-hyperconvex $T_0$-quasi-metric space and let $A\subseteq X$ be nonempty and bounded.
\begin{enumerate}[label=(\alph*)]
\item There is a unique point $\bmeet A\in X$ such that
\begin{equation}\label{eq:meet}
d\Bigl(w,\bmeet A\Bigr)=\sup_{a\in A}d(w,a)\qquad\text{for all } w\in X,
\end{equation}
and $\bmeet A$ is the infimum of $A$ in $(X,\le_d)$.
\item There is a unique point $\bjoin A\in X$ such that
\begin{equation}\label{eq:join}
d\Bigl(\bjoin A,w\Bigr)=\sup_{a\in A}d(a,w)\qquad\text{for all } w\in X,
\end{equation}
and $\bjoin A$ is the supremum of $A$ in $(X,\le_d)$.
\end{enumerate}
\end{theorem}

\begin{proof}
Let $\delta=\diam A$, fix $a_0\in A$ and put $m(w)=\sup_{a\in A}d(w,a)\le d(w,a_0)+\delta$. Apply \eqref{eq:qh} to the family indexed by $A\sqcup X$ consisting of the points $a\in A$ with radii $(\delta,0)$ and the points $w\in X$ with radii $(m(w),d(a_0,w))$. The hypothesis of \eqref{eq:qh} holds: $d(a,a')\le\delta$; $d(a,w)\le d(a,a_0)+d(a_0,w)\le\delta+d(a_0,w)$; $d(w,a)\le m(w)$; and $d(w,w')\le d(w,a_0)+d(a_0,w')\le m(w)+d(a_0,w')$. A point $z$ of the intersection satisfies $d(z,a)=0$ for $a\in A$ (so $z\le_d a$ for every $a\in A$, and in particular $\sup_{a\in A}d(z,a)=0$) and $d(w,z)\le m(w)$ for $w\in X$. Since $d(w,a)\le d(w,z)+d(z,a)=d(w,z)$, we get $d(w,z)=m(w)$, which is \eqref{eq:meet}. If $z'$ also satisfies \eqref{eq:meet}, then substituting $w=z'$ in the identity gives $\sup_{a\in A}d(z',a)=d(z',z')=0$, so in particular $d(z',z)\le\sup_ad(z',a)=0$; and substituting $w=z$ in the identity for $z'$ gives $d(z,z')=\sup_{a\in A}d(z,a)=0$ since $d(z,a)=0$ for every $a\in A$. Thus $z=z'$ by $T_0$. The point $z$ is a lower bound of $A$, and if $w$ is a lower bound, then $m(w)=0$, so $w\le_d z$. Part (b) follows by applying (a) to the $q$-hyperconvex space $(X,d^{-1})$, whose specialization order is the reverse of $\le_d$: the point $\bjoin A$ is the metric meet of $A$ in $(X,d^{-1})$, and $d^{-1}(w,\bjoin A)=\sup_ad^{-1}(w,a)$ is \eqref{eq:join}.
\end{proof}

We call $\bmeet A$ and $\bjoin A$ the \emph{metric meet} and the \emph{metric join} of $A$, and we write $x\wedge y=\bmeet\{x,y\}$ and $x\vee y=\bjoin\{x,y\}$. Formula \eqref{eq:meet} is stronger than the statement that $\bmeet A$ is an infimum: it determines all distances \emph{to} $\bmeet A$. Note also that the radii in the proof are explicit and finite, which is why boundedness of $A$, but not of $X$, is needed.

\begin{example}\label{ex:meets}
\begin{enumerate}[label=(\alph*)]
\item In $(\R,u)$, \eqref{eq:meet} reads $(w-\inf A)\dotminus0=\sup_{a\in A}(w-a)\dotminus0$, so $\bmeet A=\inf A$ and $\bjoin A=\sup A$ for every nonempty bounded $A\subseteq\R$. The set $A=\R$ has no lower bound, which shows that boundedness of $A$ cannot be dropped.
\item In $([0,1]^K,u_\infty)$ the metric meet and join are the coordinatewise infimum and supremum, because $\sup_k\sup_{a}(w_k-a_k)\dotminus0=\sup_a\sup_k(w_k-a_k)\dotminus0$.
\item In the Isbell-hull of the symmetric two-point space $\{a,b\}$ --- a square isometric to $([0,1]^2,u_\infty)$ with $a=(1,0)$ and $b=(0,1)$, computed in \cite[Ex.~5.6]{Survey} --- we have $a\wedge b=(0,0)$ and $a\vee b=(1,1)$, and $u_\infty(w,a\wedge b)=\max\{w_1,w_2\}=\max\{u_\infty(w,a),u_\infty(w,b)\}$ for all $w\in[0,1]^2$.
\end{enumerate}
\end{example}

\begin{proposition}\label{prop:powers}
Let $(X,d)$ be a $q$-hyperconvex $T_0$-quasi-metric space, $c\in X$ and $r\ge0$. There are unique points $r\pitchfork c$ and $r\ast c$ with $d(w,r\pitchfork c)=d(w,c)\dotminus r$ and $d(r\ast c,w)=d(c,w)\dotminus r$ for all $w\in X$. Moreover $c\le_d r\pitchfork c$, $d(r\pitchfork c,c)\le r$, $r\ast c\le_d c$ and $d(c,r\ast c)\le r$.
\end{proposition}

\begin{proof}
Apply \eqref{eq:qh} to the point $c$ with radii $(0,r)$ and the points $w\in X$ with radii $(d(w,c)\dotminus r,\ d(c,w)+r)$; the hypothesis holds since $d(w,w')\le d(w,c)+d(c,w')\le(d(w,c)\dotminus r)+r+d(c,w')$. A point $p$ of the intersection satisfies $d(c,p)=0$, $d(p,c)\le r$ and $d(w,p)\le d(w,c)\dotminus r$, and $d(w,c)\le d(w,p)+r$ gives equality. For uniqueness, if $p$ and $p'$ both have the defining property, then substituting $w=p'$ in the identity for $p$ gives $d(p',p)=d(p',c)\dotminus r$; and since $p'$ also satisfies the double ball condition with radii $(0,r)$ at $c$, we have $d(p',c)\le r$, so $d(p',c)\dotminus r=0$ and $d(p',p)=0$. Symmetrically $d(p,p')=0$, so $p=p'$. Finally $d(c,p)=d(c,c)\dotminus r=0$ gives $c\le_dp$. Copowers are powers in the $q$-hyperconvex space $(X,d^{-1})$, and the last two properties follow by conjugation.
\end{proof}

In $(\R,u)$ we have $r\pitchfork c=c+r$ and $r\ast c=c-r$, and in $([0,1]^K,u_\infty)$ we have $r\pitchfork c=\min\{c+r\mathbf 1,\mathbf 1\}$ and $r\ast c=\max\{c-r\mathbf 1,\mathbf 0\}$, computed coordinatewise. Copowers and powers have a simple order-theoretic description.

\begin{lemma}\label{lem:copower-least}
Let $(X,d)$ be a $q$-hyperconvex $T_0$-quasi-metric space, $c\in X$ and $r\ge0$. Then $r\ast c$ is the least element of $C_d(c,r)=\{y:d(c,y)\le r\}$, and $r\pitchfork c$ is the greatest element of $C_{d^{-1}}(c,r)=\{y:d(y,c)\le r\}$.
\end{lemma}

\begin{proof}
Taking $w=r\ast c$ in the defining identity $d(r\ast c,w)=d(c,w)\dotminus r$ gives $d(c,r\ast c)\dotminus r=0$, so $r\ast c\in C_d(c,r)$. If $d(c,y)\le r$, then $d(r\ast c,y)=d(c,y)\dotminus r=0$, i.e.\ $r\ast c\le_dy$. The statement on powers follows by passing to $(X,d^{-1})$.
\end{proof}

\begin{corollary}\label{cor:complete-lattice}
Let $(X,d)$ be a $q$-hyperconvex $T_0$-quasi-metric space. Every nonempty $A\subseteq X$ that is bounded in the quasi-metric sense $\diam A<\infty$ has an infimum and a supremum for $\le_d$. In particular, $(X,\le_d)$ is a conditionally complete lattice: every nonempty subset bounded in the quasi-metric sense has both an infimum and a supremum, which agree with the metric meet and join of \cref{thm:meets}. If $X$ itself is bounded, then $(X,\le_d)$ is a complete lattice with least element $\bot=\bmeet X$ and greatest element $\top=\bjoin X$.
\end{corollary}

\begin{remark}\label{rem:fujii}
A Lawvere metric space in the sense of \cite{Lawvere1973} is a set $X$ equipped with $d\colon X\times X\to[0,\infty]$ satisfying $d(x,x)=0$ and $d(x,z)\le d(x,y)+d(y,z)$; following \cite[Definition~3.6 and Proposition~7.6]{Fujii2020}, it is \emph{Isbell-convex} if for every family $(x_i,r_i,s_i)_{i\in I}$ with $r_i,s_i\in[0,\infty]$ and $d(x_i,x_j)\le r_i+s_j$ for all $i,j$, the intersection $\bigcap_i(C_d(x_i,r_i)\cap C_{d^{-1}}(x_i,s_i))$ is nonempty. Fujii proves that a Lawvere metric space is Isbell-convex if and only if it is cocomplete as a $[0,\infty]$-enriched category \cite[Proposition~7.6]{Fujii2020}; applied to a bounded $q$-hyperconvex $T_0$-quasi-metric space, the forward direction yields \cref{cor:complete-lattice} and \cref{prop:powers} because \eqref{eq:meet} says that $\bmeet A$ is a conical limit \cite{Stubbe2006}. What \cref{thm:meets,prop:powers} add is twofold. The proofs here use only the definition of $q$-hyperconvexity with finite radii, and the explicit radii they produce are reused in \cref{thm:tarski-qh,prop:cne-qh,prop:comparable} below. The formulas also extend to bounded subsets of unbounded $q$-hyperconvex spaces, which lie outside Fujii's framework: $(\R,u)$ with $u(x,y)=x\dotminus y$ is $q$-hyperconvex, but the family $(x,0,\infty)_{x\in\R}$ shows that it is not Isbell-convex as a Lawvere metric space.

In a related direction, Conradie, K\"unzi and Olela Otafudu describe a vector-lattice structure on the Isbell-convex hull of an asymmetrically normed real vector space \cite{CKO2017}. The present results and theirs are adjacent but not comparable. The vector-lattice structure of \cite{CKO2017} assumes an ambient linear structure and shows that the lattice operations are compatible with vector addition and scalar multiplication; \cref{thm:meets} assumes only a $q$-hyperconvex $T_0$-quasi-metric structure, with no vector operations in sight, and gives metric formulas for the lattice operations in a strictly broader class of spaces. On the overlap, namely bounded $q$-hyperconvex asymmetrically normed spaces, the two underlying lattices coincide, because both are the specialization order and both are characterized by their respective metric formulas.

The hypothesis $d\colon X\times X\to[0,\infty)$ is not needed for \cref{thm:meets}: the proof applies verbatim to an extended $T_0$-quasi-pseudometric $d\colon X\times X\to[0,\infty]$ provided \eqref{eq:qh} is still assumed for families with finite radii, and provided $\bmeet A$ and $\bjoin A$ are taken among the points of $X$ at finite distance from $A$ in both directions.
\end{remark}

\begin{corollary}\label{cor:T1}
A $T_1$-quasi-metric space with at least two points is not $q$-hyperconvex. In particular a nonempty metric space is $q$-hyperconvex if and only if it is a singleton.
\end{corollary}

\begin{proof}
Let $X$ be $q$-hyperconvex and $T_1$, and let $x\ne y$ in $X$. The set $\{x,y\}$ is bounded, so by \cref{thm:meets} the point $z=x\wedge y$ exists and satisfies $z\le_dx$ and $z\le_dy$. Since $\le_d$ is the identity relation (\cref{lem:spec}(c)), $z=x$ and $z=y$, a contradiction. A singleton is clearly $q$-hyperconvex.
\end{proof}

\begin{example}\label{ex:sorgenfrey}
On $\R$ let $d(x,y)=y-x$ if $x\le y$ and $d(x,y)=1$ if $x>y$. This is a $T_1$-quasi-metric whose conjugate topology is the Sorgenfrey topology, generated by the intervals $(x-\varepsilon,x]$. By \cref{cor:T1} it is not $q$-hyperconvex; concretely, $\{0,1\}$ has no lower bound other than points $z$ with $d(z,0)=d(z,1)=0$, and there are none. In the bounded case \cref{cor:T1} can also be derived from Fujii's description of Isbell-convex spaces as complete Lawvere metric spaces \cite{Fujii2020}, since the underlying order of a $T_1$-space is trivial.
\end{example}

The next result characterizes $q$-hyperconvexity of bounded spaces in terms of the lattice structure and of copowers and powers. For bounded spaces it is closely related to Fujii's theorem that Isbell-convexity is equivalent to completeness \cite[Proposition~7.6]{Fujii2020}; the proof below is elementary and uses only the defining identities of joins and copowers.

\begin{theorem}\label{thm:char}
For a bounded $T_0$-quasi-metric space $(X,d)$ the following are equivalent:
\begin{enumerate}[label=(\roman*)]
\item $(X,d)$ is $q$-hyperconvex;
\item $(X,\le_d)$ is a complete lattice whose infima and suprema satisfy \eqref{eq:meet} and \eqref{eq:join}, and the power $r\pitchfork c$ and the copower $r\ast c$ exist for all $c\in X$ and $r\ge0$.
\end{enumerate}
\end{theorem}

\begin{proof}
The implication (i)$\Rightarrow$(ii) is \cref{thm:meets,prop:powers}. Conversely, assume (ii), and let $(x_i,r_i,s_i)_{i\in I}$ be a nonempty family with $d(x_i,x_j)\le r_i+s_j$ for all $i$ and $j$. Put $a_i=r_i\ast x_i$, the least element of $C_d(x_i,r_i)$ (\cref{lem:copower-least}), and $z=\bjoin_ia_i$. Then $z\ge_da_i$, so $d(x_i,z)\le d(x_i,a_i)+d(a_i,z)\le r_i$. By \eqref{eq:join} and the definition of copowers,
\[
d(z,x_j)=\sup_id(a_i,x_j)=\sup_i\bigl(d(x_i,x_j)\dotminus r_i\bigr)\le s_j .
\]
Hence $z$ lies in every double ball, and $(X,d)$ is $q$-hyperconvex.
\end{proof}

The proof of (ii)$\Rightarrow$(i) only uses joins and copowers, and the dual argument only uses meets and powers. Neither part of (ii) can be omitted. The subspace $\{0,1\}$ of $([0,1],u)$ is a complete lattice whose infima and suprema satisfy \eqref{eq:meet} and \eqref{eq:join}, but the copower $\frac12\ast1$ does not exist in it, since it would be a point $p$ with $u(p,0)=u(1,0)\dotminus\frac12=\frac12$. The two-point metric space is not a lattice. Boundedness also lets us extend the power and copower formulas to $r=\infty$: $\infty\pitchfork c=\top$ and $\infty\ast c=\bot$, since $d(w,\top)\le\diam X\le\infty=d(w,c)\dotminus\infty$ (interpreting $t\dotminus\infty=0$) and dually for $\bot$.

\begin{proposition}\label{prop:meets-nonexp}
Let $(X,d)$ be $q$-hyperconvex and $T_0$, and let $A,B\subseteq X$ be nonempty and bounded. Then $d(\bmeet A,\bmeet B)\le\sup_{b\in B}\inf_{a\in A}d(a,b)$ and $d(\bjoin A,\bjoin B)\le\sup_{a\in A}\inf_{b\in B}d(a,b)$. In particular $\wedge$ and $\vee$ are nonexpansive for the maximum quasi-metric on $X\times X$.
\end{proposition}

\begin{proof}
Applying \eqref{eq:meet} with $A$ replaced by $B$ and $w=\bmeet A$, and using $d(\bmeet A,a)=0$ for $a\in A$, we obtain
\begin{align*}
d\Bigl(\bmeet A,\bmeet B\Bigr)&=\sup_{b\in B}d\Bigl(\bmeet A,b\Bigr)\le\sup_{b\in B}\inf_{a\in A}\Bigl(d\bigl(\bmeet A,a\bigr)+d(a,b)\Bigr)\\
&=\sup_{b\in B}\inf_{a\in A}d(a,b).
\end{align*}
Similarly, \eqref{eq:join} with $w=\bjoin B$ and $d(b,\bjoin B)=0$ for $b\in B$ give
\begin{align*}
d\Bigl(\bjoin A,\bjoin B\Bigr)&=\sup_{a\in A}d\Bigl(a,\bjoin B\Bigr)\le\sup_{a\in A}\inf_{b\in B}\Bigl(d(a,b)+d\bigl(b,\bjoin B\bigr)\Bigr)\\
&=\sup_{a\in A}\inf_{b\in B}d(a,b).
\end{align*}
For the binary case take $A=\{x,y\}$ and $B=\{x',y'\}$, and bound the infima by $d(x,x')$ and $d(y,y')$.
\end{proof}

The right-hand sides in \cref{prop:meets-nonexp} are the two components of the Hausdorff quasi-pseudometric on nonempty subsets used in \cite{KO2012}. Thus $\bmeet$ and $\bjoin$ are nonexpansive maps from the bounded nonempty subsets, with the appropriate component of the Hausdorff quasi-pseudometric, to $X$.

\begin{proposition}\label{prop:admissible}
Let $(X,d)$ be $q$-hyperconvex and $T_0$, and let $D=\bigcap_i(C_d(x_i,r_i)\cap C_{d^{-1}}(x_i,s_i))$ be nonempty. Then $\bmeet A,\bjoin A\in D$ for every nonempty bounded $A\subseteq D$.
\end{proposition}

\begin{proof}
Fix $a\in A$. By \eqref{eq:meet}, $d(x_i,\bmeet A)=\sup_{a'\in A}d(x_i,a')\le r_i$, and $d(\bmeet A,x_i)\le d(\bmeet A,a)+d(a,x_i)=d(a,x_i)\le s_i$. Dually, by \eqref{eq:join}, $d(\bjoin A,x_i)=\sup_{a'}d(a',x_i)\le s_i$ and $d(x_i,\bjoin A)\le d(x_i,a)+d(a,\bjoin A)\le r_i$.
\end{proof}

In particular order intervals $[a,b]={\uparrow}a\cap{\downarrow}b$ are complete sublattices, being double balls with one radius $0$. We shall see in \cref{thm:not-sublattice} that $q$-hyperconvex subspaces which are not of this form need not be sublattices.

\section{Fixed points of order-preserving maps}\label{sec:fixed}

The lattice structure of \cref{sec:lattice} combines with the fixed point theorem of Knaster and Tarski. We recall the theorem with a proof, including the version for commuting families, derive fixed points of all order-preserving self-maps of bounded $q$-hyperconvex spaces, describe the extremal fixed points metrically, compare the result with the metric approach to nonexpansive maps, and give examples.

\begin{theorem}[Knaster--Tarski \cite{Tarski1955}]\label{thm:tarski}
Let $L$ be a complete lattice.
\begin{enumerate}[label=(\alph*)]
\item If $f\colon L\to L$ is order-preserving, then $\Fix(f)$ is a nonempty complete lattice in the induced order, with least element $\inf\{x:f(x)\le x\}$ and greatest element $\sup\{x:x\le f(x)\}$.
\item If $\mathcal F$ is a family of pairwise commuting order-preserving self-maps of $L$, then the set $\Fix(\mathcal F)$ of common fixed points is a nonempty complete lattice in the induced order.
\end{enumerate}
\end{theorem}

\begin{proof}
(b) Let $P=\{x\in L:x\le f(x)\text{ for all }f\in\mathcal F\}$, which contains the least element of $L$. If $x\in P$ and $g\in\mathcal F$, then for every $f\in\mathcal F$ we have $g(x)\le g(f(x))=f(g(x))$, since $x\le f(x)$ and $g$ is order-preserving; hence $g(x)\in P$. Let $p=\sup P$. For $x\in P$ and $f\in\mathcal F$ we have $x\le f(x)\le f(p)$, so $p\le f(p)$ for all $f$, i.e.\ $p\in P$. Then $f(p)\in P$ for every $f$, so $f(p)\le p$. Thus $p$ is a common fixed point, and it is the greatest one because $\Fix(\mathcal F)\subseteq P$.

Now let $S\subseteq\Fix(\mathcal F)$ and $t=\sup S$ in $L$. If $t\le x$, then $s=f(s)\le f(x)$ for all $s\in S$ and $f\in\mathcal F$, so $t\le f(x)$. Hence every $f\in\mathcal F$ maps the complete lattice $[t,\top]$ into itself, and the maps still commute. Applying the dual of the previous paragraph to $[t,\top]$, with $P$ replaced by $\{x\in[t,\top]:f(x)\le x\text{ for all }f\}$ and $\sup$ by $\inf$, we obtain a least common fixed point $q\ge t$. Every common fixed point $y$ above all elements of $S$ lies in $[t,\top]$, so $q\le y$. Thus $q$ is the supremum of $S$ in $\Fix(\mathcal F)$. Infima in $\Fix(\mathcal F)$ are obtained dually, so $\Fix(\mathcal F)$ is a complete lattice.

(a) is the case $\mathcal F=\{f\}$ of (b); the formulas for the extremal fixed points are the sets $P$ and its dual from the proof.
\end{proof}

\begin{theorem}\label{thm:tarski-qh}
Let $(X,d)$ be a bounded $q$-hyperconvex $T_0$-quasi-metric space and let $T\colon X\to X$ be order-preserving for $\le_d$; no continuity of $T$ is assumed.
\begin{enumerate}[label=(\alph*)]
\item $\Fix(T)$ is a nonempty complete lattice in the order induced by $\le_d$.
\item The least fixed point $\ell$ and the greatest fixed point $g$ are given by $\ell=\bmeet\{x:Tx\le_dx\}$ and $g=\bjoin\{x:x\le_dTx\}$; equivalently,
\[
d(w,\ell)=\sup\{d(w,x):Tx\le_dx\},\qquad d(g,w)=\sup\{d(x,w):x\le_dTx\}\qquad(w\in X).
\]
\item If $x_0\le_dTx_0$, then $T$ has a least fixed point in ${\uparrow}x_0$; if $Tx_0\le_dx_0$, then $T$ has a greatest fixed point in ${\downarrow}x_0$.
\item Every family of pairwise commuting order-preserving self-maps of $X$ has a nonempty complete lattice of common fixed points.
\end{enumerate}
\end{theorem}

\begin{proof}
By \cref{cor:complete-lattice}, $(X,\le_d)$ is a complete lattice, and its infima and suprema are the metric meets and joins of \cref{thm:meets}. Parts (a), (b) and (d) follow from \cref{thm:tarski}; the metric formulas in (b) are \eqref{eq:meet} and \eqref{eq:join}. For (c), the set ${\uparrow}x_0$ is the interval $[x_0,\top]$ of a complete lattice, hence a complete lattice. If $x_0\le_dx$, then $x_0\le_dTx_0\le_dTx$, so $T$ maps ${\uparrow}x_0$ into itself, and (a) applies to the restriction. The second statement is dual.
\end{proof}

By \cref{lem:spec}(a), \cref{thm:tarski-qh} applies to nonexpansive maps. The following corollary compares the result with the metric approach.

\begin{corollary}\label{cor:nonexp}
Let $(X,d)$ be a bounded $q$-hyperconvex $T_0$-quasi-metric space and let $T\colon X\to X$ be nonexpansive.
\begin{enumerate}[label=(\alph*)]
\item $\Fix(T)$ is a nonempty complete lattice for $\le_d$, with least and greatest elements as in \cref{thm:tarski-qh}.
\item $\Fix(T)$ is $q$-hyperconvex \cite[Theorem~3.3]{KO2012}, and there is a nonexpansive retraction of $X$ onto $\Fix(T)$.
\item For every nonempty $A\subseteq\Fix(T)$, the infimum $\bmeet^{F}A$ of $A$ in $F=\Fix(T)$ satisfies $d(w,\bmeet^FA)=\sup_{a\in A}d(w,a)$ for all $w\in F$, and dually for suprema.
\end{enumerate}
\end{corollary}

\begin{proof}
(a) follows from \cref{thm:tarski-qh,lem:spec}(a). In (b), the $q$-hyperconvexity of $\Fix(T)$ is \cite[Theorem~3.3]{KO2012}, and the retraction exists by \cref{thm:injective}. For (c), the specialization order of the subspace $(F,d)$ is the restriction of $\le_d$, and $F$ is bounded and $q$-hyperconvex; apply \cref{thm:meets} to $F$.
\end{proof}

\begin{remark}\label{rem:compare}
\Cref{thm:tarski-qh,cor:nonexp} and \cite[Theorem~3.3]{KO2012} together give three distinct pieces of information about the fixed point set of a nonexpansive self-map $T$ of a bounded $q$-hyperconvex $T_0$-quasi-metric space $(X,d)$.
\begin{enumerate}[leftmargin=2.2em,itemsep=0.25em]
\item[$(\alpha)$] \emph{Existence of a fixed point, for every order-preserving $T$} --- including discontinuous $T$ --- is \cref{thm:tarski-qh}(a). This is new; it uses neither the metric Sine--Soardi theorem nor nonexpansiveness. The proof extracts from the metric structure the order-theoretic input that Knaster--Tarski consumes.
\item[$(\beta)$] \emph{The fixed point set is $q$-hyperconvex, for every nonexpansive $T$}, is \cite[Theorem~3.3]{KO2012}, reproduced here as \cref{cor:nonexp}(b). This is not reproven: we quote it. Its proof passes to $(X,d^s)$ and uses the Baillon/Sine--Soardi line \cite{EK2001}.
\item[$(\gamma)$] \emph{The fixed point set is a complete lattice whose infima and suprema are the metric meets and joins of \cref{thm:meets}, for every order-preserving $T$}, is \cref{thm:tarski-qh}(b) combined with \cref{cor:nonexp}(c). The complete-lattice statement goes through Knaster--Tarski and is new in the $q$-hyperconvex setting; the identification of the lattice operations as metric meets and joins in $\Fix(T)$ itself uses \cref{thm:meets} applied to the subspace $\Fix(T)$, which inherits $q$-hyperconvexity from $(\beta)$.
\end{enumerate}
Thus $(\alpha)$ replaces nonexpansiveness by order-preservation in the hypotheses, $(\beta)$ is quoted from \cite{KO2012}, and $(\gamma)$ refines the conclusion with explicit metric lattice operations in $\Fix(T)$. The proof of ($\alpha$) uses only \cref{thm:meets}, which follows directly from the definition of $q$-hyperconvexity, together with the Knaster--Tarski theorem; neither the metric theory nor the axiom of choice enters. The lattice operations of $\Fix(T)$ provided by ($\gamma$) need not coincide with those of $X$ (\cref{thm:not-sublattice}).
\end{remark}

The extremal fixed points depend monotonically on the map.

\begin{proposition}\label{prop:monotone-dependence}
Let $(X,d)$ be a bounded $q$-hyperconvex $T_0$-quasi-metric space and let $T,S\colon X\to X$ be order-preserving with $Tx\le_dSx$ for all $x\in X$. Then $\ell_T\le_d\ell_S$ and $g_T\le_dg_S$, where $\ell$ and $g$ denote least and greatest fixed points.
\end{proposition}

\begin{proof}
If $Sx\le_dx$, then $Tx\le_dSx\le_dx$; hence $\{x:Sx\le_dx\}\subseteq\{x:Tx\le_dx\}$, and taking infima gives $\ell_T\le_d\ell_S$. Similarly $\{x:x\le_dTx\}\subseteq\{x:x\le_dSx\}$ gives $g_T\le_dg_S$.
\end{proof}

\begin{proposition}\label{prop:prefixed}
Let $(X,d)$ be a bounded $q$-hyperconvex $T_0$-quasi-metric space and $T\colon X\to X$ order-preserving. Let $P=\{x:x\le_dTx\}$ and $Q=\{x:Tx\le_dx\}$.
\begin{enumerate}[label=(\alph*)]
\item $T(P)\subseteq P$, $T(Q)\subseteq Q$, $P$ is closed under arbitrary joins and $Q$ under arbitrary meets; in particular $P$ and $Q$ are complete lattices in the induced order.
\item If $T$ is nonexpansive, then $P$ and $Q$ are closed in $(X,d^s)$.
\end{enumerate}
\end{proposition}

\begin{proof}
(a) If $x\le_dTx$, then $Tx\le_dT(Tx)$, so $Tx\in P$. If $A\subseteq P$ and $s=\bjoin A$, then $a\le_dTa\le_dTs$ for all $a\in A$, so $s\le_dTs$; the empty join $\bot$ lies in $P$. The statements for $Q$ are dual.
(b) $P=\{x:d(x,Tx)=0\}$. If $T$ is nonexpansive, it is nonexpansive for $d^s$, and $x\mapsto d(x,Tx)$ is continuous on $(X,d^s)$ because $|d(x,Tx)-d(y,Ty)|\le d^s(x,y)+d^s(Tx,Ty)\le2d^s(x,y)$ by \cref{lem:spec}(b). The same argument applies to $Q$.
\end{proof}

Examples showing what the theorem does and does not give follow.

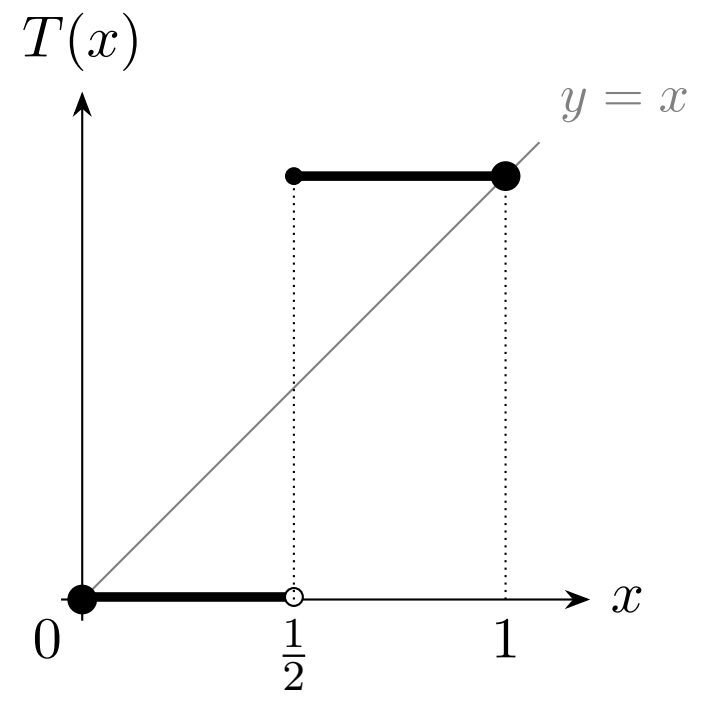
\begin{figure}[ht]
\centering
\begin{tikzpicture}[scale=2.8]
 \draw[-{Stealth}] (-0.05,0) -- (1.2,0) node[right]{$x$};
 \draw[-{Stealth}] (0,-0.05) -- (0,1.2) node[above]{$T(x)$};
 \draw[gray] (0,0) -- (1.08,1.08) node[above right,font=\small]{$y=x$};
 \draw[line width=1.8pt] (0,0.006) -- (0.5,0.006);
 \draw[fill=white] (0.5,0.006) circle (0.6pt);
 \draw[line width=1.8pt] (0.5,1) -- (1,1);
 \fill (0.5,1) circle (0.6pt);
 \fill (0,0) circle (1.0pt) node[below left]{$0$};
 \fill (1,1) circle (1.0pt);
 \draw[densely dotted] (1,0) node[below]{$1$} -- (1,1);
 \draw[densely dotted] (0.5,0) node[below]{$\tfrac12$} -- (0.5,1);
\end{tikzpicture}
\caption{The order-preserving step map of \cref{ex:step}: its fixed points $0$ and $1$ form a complete lattice which is not $q$-hyperconvex.}
\label{fig:step}
\end{figure}

\begin{example}[A discontinuous order-preserving map]\label{ex:step}
Let $X=[0,1]$ with $u$, which is bounded and $q$-hyperconvex by \cref{lem:products}, and let $T(x)=0$ for $x<\frac12$ and $T(x)=1$ for $x\ge\frac12$. Then $T$ is order-preserving and discontinuous, and $\Fix(T)=\{0,1\}$, a complete lattice with $\ell=0$ and $g=1$, as predicted by \cref{thm:tarski-qh}. However $\{0,1\}$ is not $q$-hyperconvex: the point $1$ with radii $(\frac12,1)$ and the point $0$ with radii $(1,\frac12)$ satisfy the hypothesis of \eqref{eq:qh}, since $u(1,0)=1\le\frac12+\frac12$, but the corresponding double balls in $[0,1]$ meet only at $\frac12$. Thus the conclusion of \cref{cor:nonexp}(b) genuinely requires nonexpansiveness (\cref{fig:step}).
\end{example}

\begin{example}[The hypotheses cannot be dropped]\label{ex:hypotheses}
\begin{enumerate}[label=(\alph*)]
\item \emph{Boundedness.} The map $T(x)=x+1$ on the unbounded $q$-hyperconvex space $(\R,u)$ is an isometry, hence order-preserving, and has no fixed point.
\item \emph{Monotonicity.} The map $T(x)=1$ for $x<\frac12$ and $T(x)=0$ for $x\ge\frac12$ on $([0,1],u)$ has no fixed point; it is order-reversing.
\item \emph{$q$-hyperconvexity.} On the two-point metric space $\{a,b\}$ the specialization order is trivial, so the isometry exchanging $a$ and $b$ is order-preserving; it has no fixed point.
\end{enumerate}
\end{example}

\begin{example}[Commuting maps of a square]\label{ex:commuting}
Let $X=[0,1]^2$ with $u_\infty$, let $\sigma(x_1,x_2)=(x_2,x_1)$, and for $0\le\alpha\le\beta\le1$ let $\kappa(x)=(\min\{\max\{x_1,\alpha\},\beta\},\min\{\max\{x_2,\alpha\},\beta\})$. Both maps are nonexpansive: $\sigma$ is an isometry, and $t\mapsto\min\{\max\{t,\alpha\},\beta\}$ is order-preserving and $1$-Lipschitz on $\R$, hence nonexpansive for $u$. Since the square $[\alpha,\beta]^2$ is symmetric, $\sigma\kappa=\kappa\sigma$. The common fixed points form the segment $\{(t,t):\alpha\le t\le\beta\}$, a complete lattice with least element $(\alpha,\alpha)$ and greatest element $(\beta,\beta)$, in accordance with \cref{thm:tarski-qh}(d). Here $\Fix(\kappa)=[\alpha,\beta]^2$ and $\Fix(\sigma)$, the diagonal, are complete sublattices of $X$ (\cref{prop:admissible}).
\end{example}

\section{The position of fixed point sets}\label{sec:position}

Fixed point sets of nonexpansive maps, although complete lattices in the induced order, need not be sublattices, and the specialization order itself need not be distributive. The Isbell-hull of the three-point equilateral space exhibits both phenomena.

\begin{lemma}\label{lem:Y}
Let $X=[0,1]^3$ with $q(g,h)=\max_k(h_k\dotminus g_k)$ and $Y=\{g\in X:$ the two largest coordinates of $g$ coincide$\}$. Then $(X,q)$ and $(Y,q)$ are bounded $q$-hyperconvex $T_0$-quasi-metric spaces, and $\le_q$ is the reverse coordinatewise order. Moreover $(Y,q)$ is isometric to the Isbell-hull of the space $\{a,b,c\}$ with $d(x,y)=1$ for $x\ne y$.
\end{lemma}

\begin{proof}
\emph{The space $(X,q)$.} Since $q(g,h)=u_\infty(h,g)$, the space $(X,q)$ is the conjugate of $([0,1]^3,u_\infty)$, which is $q$-hyperconvex by \cref{lem:products}; hence $(X,q)$ is $q$-hyperconvex by \cref{prop:basic}(a). Moreover $g\le_qh$ if and only if $h_k\le g_k$ for all $k$.

\emph{Double balls.} For $y\in X$ and $r,s\ge0$, the condition $q(y,h)\le r$ means $h_k\le y_k+r$ for all $k$, and $q(h,y)\le s$ means $h_k\ge y_k-s$ for all $k$. Hence the double ball of $(X,q)$ at $y$ with radii $(r,s)$ is $X\cap\prod_k[y_k-s,\,y_k+r]$.

\emph{The space $(Y,q)$.} Let $(y_i,r_i,s_i)_{i\in I}$ be a nonempty family in $Y$ with $q(y_i,y_j)\le r_i+s_j$ for all $i,j$, that is, $y_{j,k}-s_j\le y_{i,k}+r_i$ for all $i,j,k$. Put
\[
L_k=\max\Bigl\{0,\ \sup_i(y_{i,k}-s_i)\Bigr\},\qquad U_k=\min\Bigl\{1,\ \inf_i(y_{i,k}+r_i)\Bigr\}.
\]
Then $L_k\le U_k$ for each $k$, and the intersection of the double balls in $X$ is the box $\prod_k[L_k,U_k]$. We must find a point of $Y$ in this box. Relabel the three coordinates globally so that $U_1\ge U_2\ge U_3$, and put $z=(U_2,U_2,U_3)$. The two largest coordinates of $z$ coincide, so $z\in Y$, and $z_2=U_2$ and $z_3=U_3$ lie in the right intervals. It remains to show $L_1\le z_1=U_2\le U_1$. For any $y\in Y$ the largest coordinate is attained at least twice, so the smallest coordinate of $y$ is bounded above by the maximum of the other two coordinates --- equivalently, $y_k\le\max_{k'\ne k}y_{k'}$ for every $k$. Applied per index at $k=1$, this gives $y_{i,1}\le\max\{y_{i,2},y_{i,3}\}$ for every $i$, \emph{regardless} of which of $y_{i,1},y_{i,2},y_{i,3}$ is smallest under the global relabeling: the inequality holds for every coordinate of every point of $Y$. Hence
\[
y_{i,1}-s_i\le\max\{y_{i,2}-s_i,\ y_{i,3}-s_i\}\le\max\{L_2,L_3\}\le\max\{U_2,U_3\}=U_2,
\]
and $L_1\le U_2$ because also $0\le U_2$.

\emph{The Isbell-hull.} Write a function pair on $\{a,b,c\}$ as $(f,g)$ with $f=(f_a,f_b,f_c)$ and $g=(g_a,g_b,g_c)$. Ampleness means $g_x+f_y\ge1$ for $x\ne y$. The minimality equations of \cref{sec:prelim} become
\[
f_y=\max_{x\ne y}(1-g_x)\dotminus0=1-\min_{x\ne y}g_x,\qquad g_x=1-\min_{y\ne x}f_y,
\]
where we used that minimal pairs take values in $[0,1]$. Substituting the first equation into the second gives $g_x=\max_{y\ne x}\min_{z\ne y}g_z$. Order the coordinates of $g$ as $g_{\sigma1}\le g_{\sigma2}\le g_{\sigma3}$. The inner minimum $\min_{z\ne y}g_z$ equals $g_{\sigma2}$ for $y=\sigma1$ and $g_{\sigma1}$ for $y=\sigma2,\sigma3$. Hence the right-hand side equals $g_{\sigma1}$ for $x=\sigma1$, and $g_{\sigma2}$ for $x=\sigma2$ and for $x=\sigma3$. Hence minimality holds exactly when $g_{\sigma3}=g_{\sigma2}$, i.e.\ $g\in Y$, and $f$ is determined by $g$. Finally, for $g,g'\in Y$ we have $f_y-f'_y=\min_{x\ne y}g'_x-\min_{x\ne y}g_x\le\max_x(g'_x-g_x)$, so the first term in $N$ is dominated by the second, and $N((f,g),(f',g'))=\max_x(g'_x\dotminus g_x)=q(g,g')$. The embedded points are $e_a=(0,1,1)$, $e_b=(1,0,1)$ and $e_c=(1,1,0)$.
\end{proof}

\begin{theorem}\label{thm:not-sublattice}
Let $X$, $q$, $Y$ be as in \cref{lem:Y}. There is a nonexpansive $T\colon X\to X$ with $\Fix(T)=Y$, and $Y$ is not a sublattice of $(X,\le_q)$: the join of $e_a=(0,1,1)$ and $e_b=(1,0,1)$ is $(0,0,1)\notin Y$ in $X$ and $(0,0,0)$ in $Y$.
\end{theorem}

\begin{proof}
Define $\rho\colon X\to X$ by
\[
\rho(g)_k=\min\Bigl\{g_k,\ \max_{j\ne k}g_j\Bigr\},\qquad k\in\{1,2,3\}.
\]
If the coordinates of $g$ in increasing order are $a\le b\le c$, the corresponding coordinates of $\rho(g)$ are $a,b,b$. Hence $\rho(X)\subseteq Y$, and $\rho$ fixes every point of $Y$. To check nonexpansiveness, put $\varepsilon=q(g,h)$. Then $h_k\le g_k+\varepsilon$ for every $k$, and therefore
\[
\rho(h)_k\le\min\Bigl\{g_k+\varepsilon,\ \max_{j\ne k}g_j+\varepsilon\Bigr\}=\rho(g)_k+\varepsilon.
\]
Consequently $q(\rho(g),\rho(h))\le\varepsilon=q(g,h)$. Thus $T=\rho$, regarded as a self-map of $X$, is nonexpansive and has $\Fix(T)=Y$.

In $X$ the order is the reverse coordinatewise order, so the join of $e_a$ and $e_b$ is their coordinatewise minimum $(0,0,1)$, whose largest coordinate is attained only once; thus $(0,0,1)\notin Y$. In $Y$, an upper bound $h$ of $e_a$ and $e_b$ satisfies $h\le(0,1,1)$ and $h\le(1,0,1)$ coordinatewise, i.e.\ $h=(0,0,t)$ with $t\in[0,1]$, and $h\in Y$ forces $t=0$. Hence the join of $e_a$ and $e_b$ in $Y$ is $(0,0,0)$.
\end{proof}

\begin{proposition}\label{prop:M3}
The points $\bot=(1,1,1)$, $e_a=(0,1,1)$, $e_b=(1,0,1)$, $e_c=(1,1,0)$ and $\top=(0,0,0)$ of $Y$ form a sublattice of $(Y,\le_q)$ isomorphic to the diamond $M_3$. Hence the specialization order of a bounded $q$-hyperconvex space need not be distributive.
\end{proposition}

\begin{proof}
The points $\bot=(1,1,1)$ and $\top=(0,0,0)$ belong to $Y$ and are its least and greatest elements. As in the proof of \cref{thm:not-sublattice}, any two of $e_a,e_b,e_c$ have join $\top$ in $Y$. A common lower bound $h$ of $e_a$ and $e_b$ in $(Y,\le_q)$ satisfies $h\ge(1,1,1)$ coordinatewise, hence $h=\bot$, and similarly for the other pairs. The points $e_a,e_b,e_c$ are pairwise incomparable. Thus the five points form a sublattice isomorphic to $M_3$, and $e_a\wedge(e_b\vee e_c)=e_a\wedge\top=e_a\ne\bot=(e_a\wedge e_b)\vee(e_a\wedge e_c)$ shows that $Y$ is not distributive.
\end{proof}

\begin{figure}[ht]
\centering
\begin{tikzpicture}[pt/.style={circle,fill,inner sep=1.7pt}]
 \node[pt,label=below:{$\bot=(1,1,1)$}] (b) at (0,0) {};
 \node[pt,label=left:{$e_a$}] (a) at (-1.6,1.3) {};
 \node[pt,label=right:{$e_b$}] (bb) at (0,1.3) {};
 \node[pt,label=right:{$e_c$}] (c) at (1.6,1.3) {};
 \node[pt,label=above:{$\top=(0,0,0)$}] (t) at (0,2.6) {};
 \draw (b)--(a)--(t) (b)--(bb)--(t) (b)--(c)--(t);
\end{tikzpicture}
\caption{The diamond $M_3$ inside the Isbell-hull $Y$ of the equilateral three-point space.}
\label{fig:M3}
\end{figure}
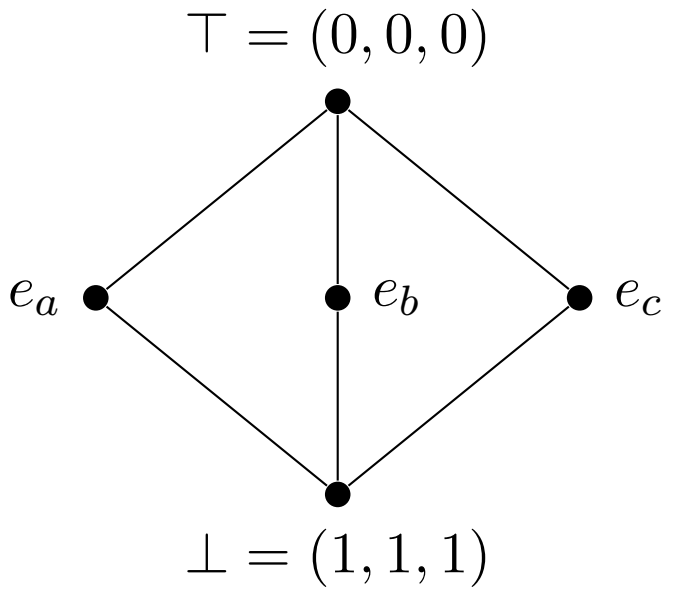

\Cref{thm:not-sublattice} should be contrasted with \cref{prop:admissible}: the fixed point set of a nonexpansive map is always $q$-hyperconvex, and a $q$-admissible set is always a complete sublattice, but $Y$ is $q$-hyperconvex without being $q$-admissible in $X$. Order-preserving self-maps of $Y$ itself behave as predicted by \cref{thm:tarski-qh}.

\begin{example}\label{ex:Y-maps}
Let $T\colon Y\to Y$, $T(g)=(m,m,m)$ with $m=\max_kg_k$. Then $T$ is order-preserving for the reverse coordinatewise order, and $q(Tg,Th)=\max_kh_k\dotminus\max_kg_k\le\max_k(h_k\dotminus g_k)=q(g,h)$, so $T$ is nonexpansive. Its fixed point set is the diagonal $\{(t,t,t):t\in[0,1]\}$, a chain from $\top$ to $\bot$ which is a complete sublattice of $Y$ and is isometric to $([0,1],u)$ via $(t,t,t)\mapsto1-t$. The least fixed point is $\bot$ and the greatest is $\top$.
\end{example}

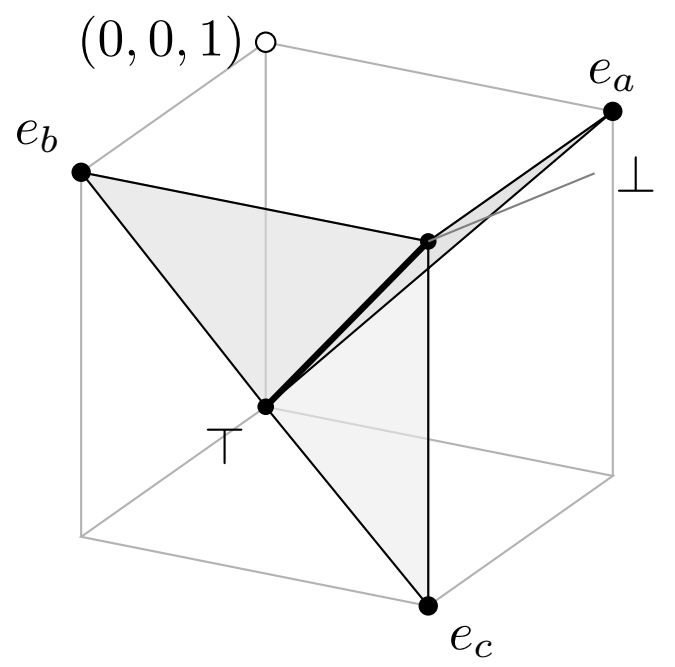
\begin{figure}[ht]
\centering
\tdplotsetmaincoords{68}{118}
\begin{tikzpicture}[tdplot_main_coords,scale=2.6]
 \draw[gray!60] (0,0,0)--(1,0,0)--(1,1,0)--(0,1,0)--cycle;
 \draw[gray!60] (0,0,1)--(1,0,1)--(1,1,1)--(0,1,1)--cycle;
 \foreach \x/\y in {0/0,1/0,1/1,0/1} \draw[gray!60] (\x,\y,0)--(\x,\y,1);
 \fill[gray!35,opacity=0.6] (0,0,0)--(1,1,1)--(0,1,1)--cycle;
 \fill[gray!25,opacity=0.6] (0,0,0)--(1,1,1)--(1,0,1)--cycle;
 \fill[gray!15,opacity=0.6] (0,0,0)--(1,1,1)--(1,1,0)--cycle;
 \draw (0,0,0)--(0,1,1)--(1,1,1) (0,0,0)--(1,0,1)--(1,1,1) (0,0,0)--(1,1,0)--(1,1,1);
 \draw[very thick] (0,0,0)--(1,1,1);
 \fill (0,0,0) circle (0.6pt) node[below left,font=\small]{$\top$};
 \fill (1,1,1) circle (0.6pt) coordinate (B);
 \node[font=\small,anchor=west,xshift=1.1cm,yshift=0.45cm] (bl) at (B) {$\bot$};
 \draw[gray,thin] (B)--(bl.west);
 \fill (0,1,1) circle (0.7pt) node[above,font=\small]{$e_a$};
 \fill (1,0,1) circle (0.7pt) node[above left,font=\small]{$e_b$};
 \fill (1,1,0) circle (0.7pt) node[below right,font=\small]{$e_c$};
 \draw[fill=white] (0,0,1) circle (0.7pt) node[left,font=\small]{$(0,0,1)$};
\end{tikzpicture}
\caption{The Isbell-hull $Y$ of the equilateral three-point space: three triangles glued along the diagonal from $\top$ to $\bot$. The join of $e_a$ and $e_b$ in the cube is $(0,0,1)\notin Y$.}
\label{fig:hull3}
\end{figure}

\section{Kleene iteration}\label{sec:kleene}

For an order-preserving self-map of a complete lattice, the least fixed point is obtained by transfinite iteration from the least element \cite{CC1979}. One asks whether $\omega$ steps suffice, i.e.\ whether $\bjoin_nT^n\bot$ is a fixed point; this holds when $T$ preserves suprema of increasing sequences.

In the theory of quasi-metric spaces as generalized domains, the relevant limits of forward-Cauchy nets are the Yoneda limits of Bonsangue, van Breugel and Rutten \cite{BvBR1998}, and nonexpansive maps which preserve them are called Yoneda-continuous \cite{BvBR1998,GL2013,GLN2017}. Nonexpansive maps need not be Yoneda-continuous: for a poset with its natural two-valued quasi-metric, the nonexpansive maps are the monotone maps, while the Yoneda-continuous ones are the Scott-continuous maps \cite[Examples~6.1 and~6.6]{GLN2017}. Such two-valued spaces are never $q$-hyperconvex unless they are singletons: if $d(y,x)=1$, the family consisting of $y$ with radii $(\frac12,1)$ and $x$ with radii $(1,\frac12)$ satisfies the hypothesis of \eqref{eq:qh}, but a point $z$ of the intersection would satisfy $d(y,z)=0=d(z,x)$, whence $d(y,x)=0$. The question is therefore what happens in the $q$-hyperconvex setting.

The main result of the section, \cref{thm:kleene}, is a five-way equivalence: for a bounded $q$-hyperconvex $T_0$-quasi-metric space, Smyth completeness, $d^s$-convergence of every increasing sequence or directed set to its supremum, preservation of increasing or directed suprema by every nonexpansive self-map, and Yoneda-continuity of every nonexpansive self-map are all equivalent, and each implies that Kleene iteration $T^n\bot$ converges to the least fixed point of every nonexpansive $T$. Two supporting results are used along the way: the supremum of a bounded directed set is its Yoneda limit (\cref{lem:orderCauchy}), and nonexpansive maps of bounded $q$-hyperconvex spaces need not preserve suprema of increasing sequences (\cref{ex:shift}).

\begin{definition}\label{def:orderCauchy}
An increasing sequence $x_0\le_d x_1\le_d\cdots$ in a $T_0$-quasi-metric space is \emph{order-Cauchy} if $\lim_n\sup_{m\ge n}d(x_m,x_n)=0$.
\end{definition}

\begin{lemma}\label{lem:orderCauchy-ds}
Let $(x_n)$ be an increasing sequence in a $T_0$-quasi-metric space $(X,d)$. Then $(x_n)$ is order-Cauchy if and only if it is $d^s$-Cauchy. More precisely, for $n\le m$,
\[
d^s(x_n,x_m)=\max\{d(x_n,x_m),d(x_m,x_n)\}=d(x_m,x_n),
\]
so $\sup_{m\ge n}d^s(x_m,x_n)=\sup_{m\ge n}d(x_m,x_n)$, and the two Cauchy conditions coincide.
\end{lemma}

\begin{proof}
For $n\le m$, $d(x_n,x_m)=0$ by the definition of $\le_d$, so $d^s(x_n,x_m)=d(x_m,x_n)$. Taking the supremum over $m\ge n$ and then the limit in $n$ establishes the equivalence.
\end{proof}

Part (b) of the following lemma, for increasing sequences, states that a $d^s$-limit of an increasing sequence is its supremum; this holds in every $T_0$-quasi-metric space and is due to L\'opez-Ram\'irez and Valero \cite{LRV2018}; see also \cite[Lemma~2.18]{EMV2019} and \cite[Proposition~2.1]{SV2018}. Parts (a) and (c) use the metric join formula \eqref{eq:join}.

\begin{lemma}\label{lem:orderCauchy}
Let $(X,d)$ be a $q$-hyperconvex $T_0$-quasi-metric space and let $D\subseteq X$ be a bounded directed set with supremum $s=\bjoin D$.
\begin{enumerate}[label=(\alph*)]
\item For $x\in D$ we have $d(x,s)=0$ and $d^s(x,s)=d(s,x)=\sup_{y\in D}d(y,x)$, and $x\mapsto d(s,x)$ is nonincreasing on $D$. Hence $D$ converges to $s$ in $(X,d^s)$ if and only if $\inf_{x\in D}d(s,x)=0$.
\item If $D$ converges in $(X,d^s)$ to a point $y$, then $y=s$.
\item $s$ is the unique Yoneda limit of the left $K$-Cauchy net $D$.
\end{enumerate}
In particular, an increasing sequence $(x_n)$ is order-Cauchy if and only if it converges to its supremum in $(X,d^s)$.
\end{lemma}

\begin{proof}
(a) By \eqref{eq:join}, $d(s,x)=\sup_{y\in D}d(y,x)$, and $d(x,s)=0$ since $s$ is an upper bound. If $x\le_d x'$ in $D$, then $d(s,x')\le d(s,x)+d(x,x')=d(s,x)$. A nonincreasing net of nonnegative reals converges to its infimum.
(b) For $x\in D$, the net lies eventually in the $d^s$-closed set ${\uparrow}x$ (\cref{lem:spec}(b)), so $x\le_d y$; and if $b$ is an upper bound of $D$, the net lies in the $d^s$-closed set ${\downarrow}b$, so $y\le_d b$. Thus $y$ is the supremum.
(c) $D$ is left $K$-Cauchy because $d(x,x')=0$ for $x\le_d x'$. For fixed $w$, $x\mapsto d(x,w)$ is nondecreasing on $D$, since $d(x,w)\le d(x,x')+d(x',w)$. Hence $\inf_{x\in D}\sup_{x'\ge x}d(x',w)=\sup_{x\in D}d(x,w)=d(s,w)$ by \eqref{eq:join}. Uniqueness follows from $T_0$.
For an increasing sequence, (a) gives $d(s,x_n)=\sup_{m\ge n}d(x_m,x_n)$, which is the order-Cauchy condition.
\end{proof}

\begin{figure}[ht]
\centering
\begin{tikzpicture}[xscale=0.42,yscale=0.9]
 \foreach \n/\y in {0/0,1/1.3,2/2.6,3/3.9}{
   \node[left,font=\small] at (-0.4,\y+0.3) {$T^{\n}\bot$};
   \foreach \k in {0,...,11}{
     \pgfmathsetmacro{\h}{(\k<\n)?0.6:0.04}
     \fill[gray!60] (\k,\y) rectangle (\k+0.7,\y+\h);
   }
   \draw (12.6,\y) rectangle (13.3,\y+0.04);
 }
 \node at (6,5.0) {$\vdots$};
 \foreach \k in {0,...,11}{\fill[gray!60] (\k,5.6) rectangle (\k+0.7,6.2);}
 \draw (12.6,5.6) rectangle (13.3,5.64);
 \node[left,font=\small] at (-0.4,5.9) {$s=\bigvee_nT^n\bot$};
 \foreach \k in {0,...,11}{\fill[gray!60] (\k,6.9) rectangle (\k+0.7,7.5);}
 \fill[black] (12.6,6.9) rectangle (13.3,7.5);
 \node[left,font=\small] at (-0.4,7.2) {$T(s)$};
 \node[font=\small] at (12.95,-0.6) {$z$};
 \node[font=\small] at (5.5,-0.6) {coordinates $x_0,x_1,x_2,\dots$};
\end{tikzpicture}
\caption{The orbit of $\bot$ in \cref{ex:shift}. Each bar is a coordinate of $x\in[0,1]^{\N}$, the last box is $z$. The orbit fills the coordinates one by one; its supremum $s$ has $z=0$, but $T(s)$ has $z=1$.}
\label{fig:shift}
\end{figure}
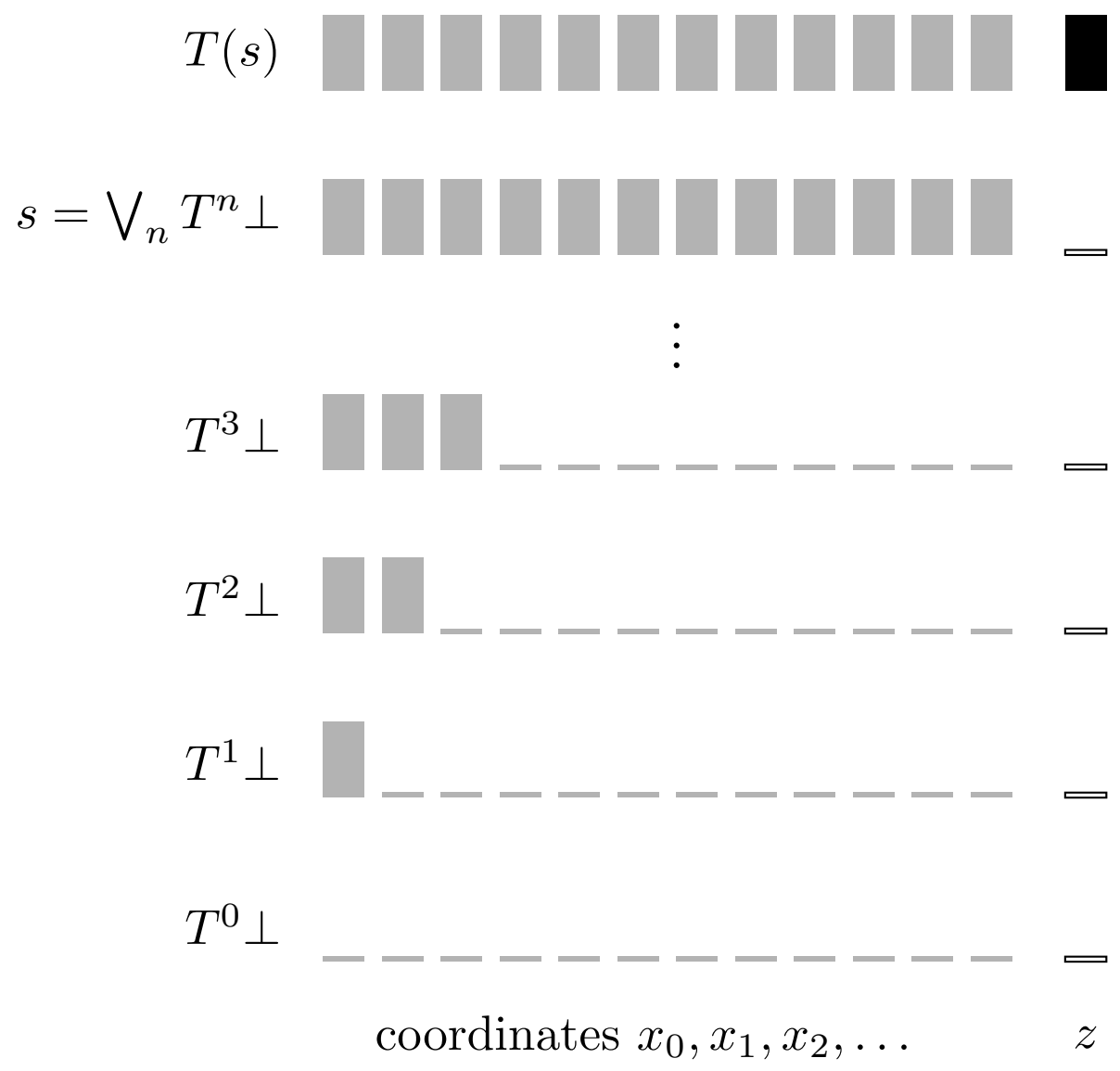

\begin{example}[Kleene iteration may need $\omega+1$ steps]\label{ex:shift}
Let $X=[0,1]^{\N}\times[0,1]$ with $d((x,z),(y,w))=\max\{\sup_k(x_k\dotminus y_k),\ z\dotminus w\}$. By \cref{lem:products}, $X$ is a bounded $q$-hyperconvex $T_0$-quasi-metric space with the coordinatewise order; its least element is $\bot=(0,0)$. Let $S(x)=(1,x_0,x_1,\dots)$ be the shift that inserts $1$, let $\Lambda(x)=\limsup_kx_k$, and define
\[
T(x,z)=\bigl(Sx,\ \max\{z,\Lambda(x)\}\bigr).
\]
\emph{$T$ is nonexpansive.} We have $\sup_k\bigl((Sx)_k\dotminus(Sy)_k\bigr)=\sup_k(x_k\dotminus y_k)$, because $(Sx)_0=(Sy)_0=1$, and $\Lambda(x)-\Lambda(y)\le\limsup_k(x_k-y_k)\le\sup_k(x_k\dotminus y_k)$. Since $\max\{a,b\}\dotminus\max\{c,e\}\le\max\{a\dotminus c,\ b\dotminus e\}$, we get $d(T(x,z),T(y,w))\le\max\{\sup_k(x_k\dotminus y_k),\ z\dotminus w\}$.

\emph{The orbit of $\bot$.} By induction $T^n\bot=(c_n,0)$, where $c_n=(1,\dots,1,0,0,\dots)$ has $n$ entries $1$, because $\Lambda(c_n)=0$. The orbit is increasing with supremum $s=(\mathbf 1,0)$, where $\mathbf 1=(1,1,\dots)$, and $d(s,T^n\bot)=1$ for all $n$.

\emph{Conclusion.} $T(s)=(\mathbf 1,\Lambda(\mathbf 1))=(\mathbf 1,1)\ne s$, so $\bjoin_nT^n\bot$ is not a fixed point and $T$ does not preserve the supremum of the increasing sequence $(T^n\bot)$. The unique fixed point is $(\mathbf 1,1)$: $Sx=x$ forces $x=\mathbf 1$, and then $z=\max\{z,1\}=1$. It equals $T(s)=T^{\omega+1}\bot$ in the transfinite iteration (\cref{fig:shift}).
\end{example}

The following theorem shows that this phenomenon is governed by Smyth completeness.

\begin{theorem}\label{thm:kleene}
Let $(X,d)$ be a bounded $q$-hyperconvex $T_0$-quasi-metric space. The following are equivalent:
\begin{enumerate}[label=(\roman*)]
\item every increasing sequence in $X$ is order-Cauchy, i.e.\ converges in $(X,d^s)$ to its supremum;
\item $(X,d)$ is Smyth complete;
\item every directed subset of $X$ converges in $(X,d^s)$ to its supremum;
\item every nonexpansive self-map of $X$ preserves suprema of directed subsets;
\item every nonexpansive self-map of $X$ preserves suprema of increasing sequences;
\item every nonexpansive self-map of $X$ is Yoneda-continuous.
\end{enumerate}
If these conditions hold, then for every nonexpansive $T\colon X\to X$ the sequence $T^n\bot$ converges in $(X,d^s)$ to the least fixed point of $T$.
\end{theorem}

\begin{proof}
We prove (i)$\Rightarrow$(ii)$\Rightarrow$(iii)$\Rightarrow$(iv)$\Rightarrow$(v)$\Rightarrow$(i) and (ii)$\Rightarrow$(vi)$\Rightarrow$(v).

(i)$\Rightarrow$(ii). By \cite[Proposition~1]{RT2015} it suffices to show that every left $K$-Cauchy sequence $(x_n)$ converges in $(X,d^s)$. Choose $n_1<n_2<\cdots$ such that $d(x_p,x_q)<2^{-k}$ whenever $n_k\le p\le q$, and put $y_k=x_{n_k}$, $\eta_k=2^{1-k}$ and $q_k=\eta_k\ast y_k$, the copower of \cref{prop:powers}. Then $q_k\le_d y_k$ and $d(y_k,q_k)\le\eta_k$, so $d^s(q_k,y_k)\le\eta_k$. The sequence $(q_k)$ is increasing:
\[
d(q_k,q_{k+1})=d(y_k,q_{k+1})\dotminus\eta_k\le\bigl(d(y_k,y_{k+1})+d(y_{k+1},q_{k+1})\bigr)\dotminus\eta_k\le(2^{-k}+2^{-k})\dotminus2^{1-k}=0 .
\]
By (i), $q_k\to s$ in $(X,d^s)$ for some $s$, hence $y_k\to s$ in $(X,d^s)$. Given $\varepsilon>0$, choose $k$ with $2^{-k}<\varepsilon$ and $d^s(y_j,s)<\varepsilon$ for all $j\ge k$. For $n\ge n_k$ we have $d(s,x_n)\le d(s,y_k)+d(y_k,x_n)<2\varepsilon$, and choosing $j\ge k$ with $n_j\ge n$ we get $d(x_n,s)\le d(x_n,y_j)+d(y_j,s)<2\varepsilon$. Thus $x_n\to s$ in $(X,d^s)$.

(ii)$\Rightarrow$(iii). A directed set $D$ is a left $K$-Cauchy net (\cref{lem:orderCauchy}(c)); by (ii) it converges in $(X,d^s)$, and its limit is $\bjoin D$ by \cref{lem:orderCauchy}(b).

(iii)$\Rightarrow$(iv). Let $T$ be nonexpansive and $D$ directed with supremum $s$. By (iii), $D$ converges to $s$ in $(X,d^s)$. Since $T$ is nonexpansive for $d^s$, the net $(Tx)_{x\in D}$ converges to $Ts$ in $(X,d^s)$. The image $T(D)$ is directed because $T$ is order-preserving (\cref{lem:spec}(a)). Applied to the net $(Tx)_{x\in D}$: for each fixed $x_0\in D$ the net lies eventually in the $d^s$-closed set ${\uparrow}Tx_0$ (\cref{lem:spec}(b)), so $Tx_0\le_d Ts$; and if $b$ is an upper bound of $T(D)$, the net lies in the $d^s$-closed set ${\downarrow}b$, so $Ts\le_d b$. Thus $Ts=\bjoin T(D)$.

(iv)$\Rightarrow$(v) is clear, since an increasing sequence is a directed set.

(v)$\Rightarrow$(i). Suppose that $(x_n)$ is increasing but not order-Cauchy, and let $s=\bjoin_nx_n$. By \cref{lem:orderCauchy}(a), $d(s,x_n)$ decreases to some $\delta>0$; in particular $s\ne x_n$ for all $n$. Let $z=\delta\ast s$, so that $z\le_d s$ and $d(s,z)\le\delta$. Moreover $z\ne s$: otherwise $d(s,w)\dotminus\delta=d(s,w)$ for all $w$, hence $d(s,w)=0$ for all $w$, contradicting $d(s,x_0)\ge\delta>0$. On $A=\{x_n:n\in\N\}\cup\{s\}$ define $f(x_n)=z$ and $f(s)=s$. Then $f$ is nonexpansive: $d(f(x_n),f(x_m))=0$, $d(f(x_n),f(s))=d(z,s)=0$, and $d(f(s),f(x_n))=d(s,z)\le\delta\le d(s,x_n)$. By \cref{thm:injective}, $f$ extends to a nonexpansive $T\colon X\to X$, and $\bjoin_nTx_n=z\ne s=T(\bjoin_nx_n)$, contradicting (v).

(ii)$\Rightarrow$(vi). Let $(x_i)$ be a left $K$-Cauchy net with Yoneda limit $y$. By (ii), $(x_i)$ converges in $(X,d^s)$ to some $y'$, which is a Yoneda limit of $(x_i)$; hence $y'=y$. Since $T$ is nonexpansive for $d^s$, $Tx_i\to Ty$ in $(X,d^s)$, so $Ty$ is the Yoneda limit of $(Tx_i)$.

(vi)$\Rightarrow$(v). By \cref{lem:orderCauchy}(c), the supremum $s$ of an increasing sequence $(x_n)$ is its Yoneda limit. By (vi), $Ts$ is the Yoneda limit of the increasing sequence $(Tx_n)$, which by \cref{lem:orderCauchy}(c) is $\bjoin_nTx_n$.

Finally, under (v) and (i), the supremum $s$ of the increasing sequence $T^n\bot$ satisfies $Ts=\bjoin_nT^{n+1}\bot=s$ and $T^n\bot\to s$ in $(X,d^s)$. If $p$ is a fixed point, then $T^n\bot\le_d p$ for all $n$ by induction, so $s\le_d p$. Thus $s$ is the least fixed point.
\end{proof}

\begin{remark}\label{rem:valero}
The implication (ii)$\Rightarrow$(i) of \cref{thm:kleene} holds in every $T_0$-quasi-metric space and is known: Smyth completeness implies that increasing sequences converge in $(X,d^s)$, and in particular that they converge in the weaker sense of the order-completeness used by Shahzad and Valero \cite[Corollary~2.16]{SV2018}; see also \cite[Proposition~2.19]{EMV2019}. In their Kleene-type theorems \cite[Theorem~2.3]{SV2018}, \cite[Corollary~2.17]{EMV2019}, order-continuity of the map along the orbit is a hypothesis. \Cref{thm:kleene} shows that for nonexpansive maps of bounded $q$-hyperconvex spaces this hypothesis is automatically satisfied if and only if the space is Smyth complete. The directions (iii)$\Rightarrow$(iv)$\Rightarrow$(v) and (ii)$\Rightarrow$(vi) are routine once (ii) is available; the content specific to the $q$-hyperconvex setting lies in the converses (i)$\Rightarrow$(ii) and (v)$\Rightarrow$(i), which use the copowers of \cref{prop:powers} and the injectivity of \cref{thm:injective}, respectively, and have no counterpart in the general quasi-metric theory of \cite{SV2018,EMV2019}. The references \cite{SV2018,EMV2019} use the conjugate convention $d_u(x,y)=\max\{y-x,0\}$ for the standard quasi-metric on $\R$.
\end{remark}

\begin{remark}\label{rem:smyth}
For general $T_0$-quasi-metric spaces, Smyth completeness is much stronger than the convergence of increasing sequences; the natural $T_0$-quasi-metric of a poset, for instance, is Smyth complete if and only if the poset satisfies the ascending chain condition, since $d^s$ is then the discrete metric. The equivalence (i)$\Leftrightarrow$(ii) uses the copowers available in $q$-hyperconvex spaces: they replace a left $K$-Cauchy sequence by an increasing sequence at small $d^s$-distance. The implication (v)$\Rightarrow$(i) uses injectivity. The map constructed there collapses a non-order-Cauchy sequence onto a copower of its supremum; in \cref{ex:shift} the failure is visible on the orbit of $\bot$ itself. Whether this is always possible is discussed in \cref{sec:conclusion}.
\end{remark}

\begin{proposition}\label{prop:compact}
Let $(X,d)$ be a bounded $q$-hyperconvex $T_0$-quasi-metric space such that $(X,d^s)$ is compact. Then every increasing sequence in $X$ is order-Cauchy; hence $X$ is Smyth complete and every nonexpansive self-map of $X$ is Yoneda-continuous. This applies in particular to boxes $\prod_{k=1}^n[a_k,b_k]$ in $(\R^n,u_\infty)$ and to the Isbell-hull of every finite $T_0$-quasi-metric space.
\end{proposition}

\begin{proof}
Let $(x_n)$ be increasing with supremum $s$. By compactness a subsequence $(x_{n_k})$ converges in $(X,d^s)$ to some point $y$. The subsequence is increasing, so $y$ is its supremum by \cref{lem:orderCauchy}(b), and this supremum is $s$ because the subsequence is cofinal. By \cref{lem:orderCauchy}(a), $d(s,x_n)$ is nonincreasing in $n$, and $d(s,x_{n_k})\le d^s(s,x_{n_k})\to0$; hence $d(s,x_n)\to0$, i.e.\ $(x_n)$ is order-Cauchy. The remaining assertions follow from \cref{thm:kleene}. For a box, $u_\infty^s$ is the maximum metric. For a finite space $X$ with $D=\diam X$, the minimal pairs satisfy $0\le f,g\le D$ and are characterized by the closed conditions of \cref{sec:prelim}; thus $Q(X)$ is a closed bounded subset of $\R^{2|X|}$, and $N^s$ is the maximum metric on pairs.
\end{proof}

In \cref{ex:shift} the space $[0,1]^{\N}\times[0,1]$ is not Smyth complete, in accordance with \cref{thm:kleene}, although each factor $[0,1]$ is.

\section{Maps that are nonexpansive on comparable pairs}\label{sec:comparable}

Monotone nonexpansive maps on metric spaces endowed with a partial order have been studied by Bachar, Bin Dehaish, Khamsi, Wi\'snicki and others \cite{BDK2016,QW2024}. In our setting the order is read off the quasi-metric, and a self-map $T$ of a $T_0$-quasi-metric space is \emph{comparably nonexpansive} if $d(Tx,Ty)\le d(x,y)$ whenever $x$ and $y$ are $\le_d$-comparable. If $x\le_d y$, the condition on the pair $(x,y)$ reads $Tx\le_d Ty$, so comparably nonexpansive maps are order-preserving; this observation, for contractions on comparable pairs, is \cite[Corollary~2.19]{SV2018}. Metric versions with an external order go back to Ran and Reurings and to Nieto and Rodr\'iguez-L\'opez \cite{RR2004,NRL2005}. In $q$-hyperconvex spaces, however, this weakening gains nothing, as the next proposition shows.

\begin{proposition}\label{prop:cne-qh}
Every comparably nonexpansive self-map of a $q$-hyperconvex $T_0$-quasi-metric space is nonexpansive.
\end{proposition}

\begin{proof}
Let $x,y\in X$ and $j=x\vee y$ (\cref{thm:meets}). Then $Tx\le_d Tj$, the pair $(j,y)$ is comparable, and $d(j,y)=\max\{d(x,y),d(y,y)\}=d(x,y)$ by \eqref{eq:join}. Hence $d(Tx,Ty)\le d(Tx,Tj)+d(Tj,Ty)\le0+d(j,y)=d(x,y)$.
\end{proof}

\begin{example}\label{ex:three}
Outside the $q$-hyperconvex setting the two notions differ. Let $X=\{\bot,a,b\}$ with $d(\bot,a)=d(\bot,b)=0$, $d(a,\bot)=d(b,\bot)=2$, $d(a,b)=1$, $d(b,a)=2$. This is a $T_0$-quasi-metric, and the map fixing $\bot$ and exchanging $a$ and $b$ is comparably nonexpansive but not nonexpansive, since $d(Ta,Tb)=2>1=d(a,b)$. In accordance with \cref{prop:cne-qh}, $\{a,b\}$ has no upper bound in $X$.
\end{example}

\begin{figure}[ht]
\centering
\begin{tikzpicture}[pt/.style={circle,fill,inner sep=1.7pt}]
 \node[pt,label=below:{$\bot$}] (b) at (0,0) {};
 \node[pt,label=left:{$a$}] (a) at (-1.8,1.7) {};
 \node[pt,label=right:{$b$}] (c) at (1.8,1.7) {};
 \draw (b)--(a) (b)--(c);
 \draw[-{Stealth}] (a) to[bend left=15] node[above,font=\small]{$d(a,b)=1$} (c);
 \draw[-{Stealth}] (c) to[bend left=15] node[below,font=\small]{$d(b,a)=2$} (a);
\end{tikzpicture}
\caption{The space of \cref{ex:three}. The swap of $a$ and $b$ is nonexpansive on the comparable pairs $\{\bot,a\}$ and $\{\bot,b\}$ but expands the pair $(a,b)$.}
\label{fig:three}
\end{figure}
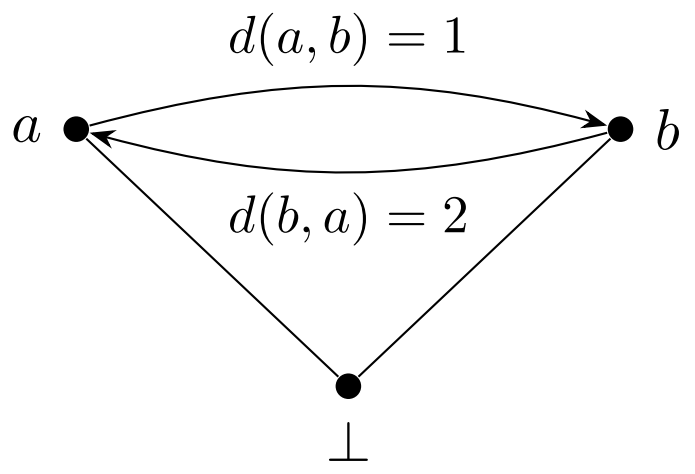

Outside $q$-hyperconvex spaces, comparably nonexpansive maps admit the following fixed point result. It is a variant of Rutten's theorem \cite{Rutten1996} and of the Kleene-type theorem of Shahzad and Valero \cite[Theorem~2.3 and Corollary~2.20]{SV2018}, in which order-continuity is replaced by the order-Cauchy property of the orbit; we include the short proof.

\begin{proposition}\label{prop:comparable}
Let $(X,d)$ be a bicomplete $T_0$-quasi-metric space and let $T\colon X\to X$ be comparably nonexpansive. Suppose that $x_0\le_dTx_0$ and that the orbit $x_n=T^nx_0$ is order-Cauchy. Then $(x_n)$ converges in $(X,d^s)$ to a point $x^*$ which is the supremum of $\{x_n:n\in\N\}$ and the least fixed point of $T$ in ${\uparrow}x_0$.
\end{proposition}

\begin{proof}
Since $T$ is order-preserving, $x_n\le_dx_{n+1}$ for all $n$ by induction. By \cref{lem:orderCauchy-ds}, the order-Cauchy sequence $(x_n)$ is $d^s$-Cauchy, and converges in the complete metric space $(X,d^s)$ to some $x^*$. For fixed $n$, all $x_m$ with $m\ge n$ lie in the $d^s$-closed set ${\uparrow}x_n$, so $x_n\le_dx^*$; and if $b$ is an upper bound of the orbit, all $x_m$ lie in the $d^s$-closed set ${\downarrow}b$, so $x^*\le_db$. Thus $x^*$ is the supremum of the orbit. From $x_n\le_dx^*$ we get $x_{n+1}\le_dTx^*$ for all $n$, and $x_0\le_dx_1$, so $x^*\le_dTx^*$. On the other hand $x_n$ and $x^*$ are comparable, so
\[
d(Tx^*,x^*)\le d(Tx^*,Tx_n)+d(x_{n+1},x^*)\le d(x^*,x_n)+0\longrightarrow0 .
\]
Hence $Tx^*\le_dx^*$, and $Tx^*=x^*$. If $p\in{\uparrow}x_0$ is a fixed point, then $x_n\le_dp$ implies $x_{n+1}=Tx_n\le_dTp=p$, so all $x_n\le_dp$ and $x^*\le_dp$.
\end{proof}

The order-Cauchy hypothesis cannot be dropped, even for nonexpansive maps of bounded $q$-hyperconvex spaces (\cref{ex:shift}), and comparable nonexpansiveness cannot be weakened to monotonicity: on $(\R,u)$ the map with $T(x)=\tfrac{x+1}{2}$ for $x<1$ and $T(x)=x+1$ for $x\ge1$ is order-preserving, the orbit of $0$ is $1-2^{-n}$, which is increasing and order-Cauchy, and $T$ has no fixed point.

\section{Open problems}\label{sec:conclusion}

Several questions are left open by the results above; for each we state what partial information or conjecture we have.

\subsection*{Kleene iteration and Scott-continuity} \Cref{thm:kleene} characterizes Smyth completeness as the equivalent of Scott-continuity for all nonexpansive maps, but the diagnostic in the proof of (v)$\Rightarrow$(i) is constructed via injectivity: it collapses an arbitrary non-order-Cauchy increasing sequence onto a copower of its supremum. In \cref{ex:shift} the failure is instead concrete on the orbit $(T^n\bot)$ of the least element.

\begin{question}\label{q:kleene}
Let $(X,d)$ be a bounded $q$-hyperconvex $T_0$-quasi-metric space that is not Smyth complete. Does there exist a nonexpansive self-map $T\colon X\to X$ whose orbit $(T^n\bot)$ is not order-Cauchy, i.e.\ such that $\bjoin_nT^n\bot$ is not a fixed point?
\end{question}

We expect the answer to be negative in general. The map constructed in the proof of (v)$\Rightarrow$(i) collapses a specified non-order-Cauchy increasing sequence onto a copower of its supremum, but there is no reason why such a sequence should arise as the orbit $(T^n\bot)$ of a nonexpansive self-map, let alone with $\bot$ as starting point. A positive answer would locate the failure of Smyth completeness on the orbit of a single canonical point.

A related question concerns the ordinals at which the transfinite Kleene iteration terminates. \Cref{ex:shift} shows that $\omega+1$ is realized in a bounded $q$-hyperconvex space; we do not know whether larger countable ordinals occur. Cousot and Cousot \cite{CC1979} provide the general transfinite framework.

\begin{question}\label{q:ordinals}
For which countable ordinals $\alpha$ does there exist a bounded $q$-hyperconvex $T_0$-quasi-metric space $(X,d)$ and a nonexpansive $T\colon X\to X$ such that the transfinite Kleene iteration from $\bot$ reaches a fixed point exactly at stage $\alpha$?
\end{question}

For unbounded $q$-hyperconvex spaces, $(\R,u)$ is not Smyth complete: the sequence $(n)_{n\in\N}$ is left $K$-Cauchy (indeed increasing) with no $u^s$-limit. Its bounded increasing sequences do converge, so the appropriate analogue of \cref{thm:kleene} presumably restricts (i), (ii), (iii) to bounded sequences and sets.

\begin{question}\label{q:unbounded}
Formulate and prove the analogue of \cref{thm:kleene} for unbounded $q$-hyperconvex $T_0$-quasi-metric spaces, restricting (i), (ii), (iii) to bounded sequences and sets.
\end{question}

\subsection*{The order-theoretic side} \Cref{prop:M3} shows that non-distributive lattices occur as specialization orders of bounded $q$-hyperconvex spaces: the diamond $M_3$ embeds in the Isbell-hull of the equilateral three-point space as a sublattice. We did not investigate modularity.

\begin{question}\label{q:modular}
Is the specialization order $(Y,\le_q)$ of the Isbell-hull of the equilateral three-point space, described in \cref{lem:Y}, modular? Does $N_5$ occur as a sublattice of $(Y,\le_q)$, or more generally of the specialization order of some bounded $q$-hyperconvex space?
\end{question}

Computational sampling of triples in $Y$ consistent with the modular law $a\le c\Rightarrow a\vee(b\wedge c)=(a\vee b)\wedge c$ returns no obvious counterexample; we conjecture that $(Y,\le_q)$ is modular but not distributive. The Banaschewski--Bruns description of the Dedekind--MacNeille completion \cite{BB1967}, combined with the correspondence of Agyingi, Haihambo and K\"unzi \cite{AHK2013} between the Dedekind--MacNeille completion of a poset and the $q$-hyperconvex hull of its natural $T_0$-quasi-metric, suggests that modularity can be read off properties of the generating poset.

\begin{question}\label{q:which-lattices}
Which complete lattices occur as the specialization order of a bounded $q$-hyperconvex $T_0$-quasi-metric space, up to isomorphism of lattices?
\end{question}

\subsection*{Sublattices among fixed point sets} \Cref{prop:admissible,thm:not-sublattice} together show that fixed point sets of nonexpansive maps are $q$-hyperconvex but need not be sublattices of the ambient space, while $q$-admissible sets are always complete sublattices. This raises:

\begin{question}\label{q:sublattice-charact}
Characterize the nonexpansive self-maps of a bounded $q$-hyperconvex $T_0$-quasi-metric space whose fixed point set is a sublattice. Equivalently, characterize the $q$-hyperconvex subspaces of a bounded $q$-hyperconvex space that are complete sublattices.
\end{question}

A sufficient condition is $q$-admissibility of the fixed point set in the sense of \cref{prop:admissible} (intersection of double balls), which holds for instance when $T$ is a composition of radial retractions. The converse is open.

\subsection*{Asymmetrically normed spaces} The $q$-hyperconvex asymmetrically normed real vector spaces we are aware of are $(\R,u)$ and $(\R^n,u_\infty)$ (and their products), together with their Isbell-convex hulls as described by Conradie, K\"unzi and Olela Otafudu \cite{CKO2017}. We do not know of further nontrivial examples.

\begin{question}\label{q:asymmetric-normed}
Are there $q$-hyperconvex asymmetrically normed real vector spaces other than $(\R^n,u_\infty)$ and their closed subspaces? More specifically, is every $q$-hyperconvex asymmetrically normed space of this form up to isometric isomorphism of normed spaces?
\end{question}

\subsection*{Applications and further directions} Valero and his coauthors used the specialization order of the complexity space to analyse recurrence equations of algorithms \cite{EMV2019}; one may ask whether the Isbell-hull of the complexity space, where \cref{thm:tarski-qh} applies without continuity assumptions, adds anything there. Analogues for $q$-spherically complete $T_0$-ultra-quasi-metric spaces, where sums of radii are replaced by maxima \cite{KO2013}, should also be within reach.

\subsection*{Acknowledgements}
[To be completed.]

\end{document}